\RequirePackage{fix-cm}
\documentclass[smallextended]{svjour3}       
\smartqed  
\usepackage{graphicx}
\usepackage{mathptmx}      

\usepackage{amssymb,amsmath,color}

\usepackage{enumitem}
\usepackage{graphicx,tikz}
\usepackage{mathtools}
\usepackage{hyperref}
\usepackage[english]{babel}
\usepackage{graphicx,epstopdf,float,tikz}
\usetikzlibrary{arrows,positioning,calc}

\usepackage{graphicx} 
\usepackage{subcaption}

\newtheorem{assumption}{Assumption}

\newcommand{\R}{\mathbb R}

\newcommand{\N}{\mathbb N}

\newcommand{\cA}{\mathcal A}

\newcommand{\cI}{\mathcal I}
\newcommand{\cJ}{\mathcal J}
\newcommand{\cK}{\mathcal K}

\newcommand{\cN}{\mathcal N}

\newcommand{\x}{\times}
\newcommand{\st}{\,:\,}

\newcommand{\inv}{^{-1}}

\DeclareMathOperator{\linspan}{span}

\DeclareMathOperator*{\subj}{\textnormal{subj.~to}}

\newcommand{\sr}{^\star} 
\newcommand{\T}{^\top}
\newcommand{\eq}{^{\rm e}}
\newcommand{\m}{_{\rm m}}
\newcommand{\p}{_{\perp}}
 
\newcommand{\tx}{\tilde{x}}
\newcommand{\tz}{\tilde{z}}
\newcommand{\txi}{\tilde{\xi}}
\newcommand{\tlam}{\tilde{\lambda}}

\allowdisplaybreaks

\newcommand{\1}{\boldsymbol{1}}

\newcommand{\sprod}[2]{\langle #1, #2\rangle}
\newcommand{\nablaL}{\Phi}

\journalname{journalname}
\begin{document}

\title{An exponentially stable discrete-time primal-dual algorithm for distributed constrained optimization\thanks{This work has been partially supported by the European Union’s Horizon 2020 Research and Innovation Programme under grant agreement No. 739551 (KIOS CoE) and by the European Union under NextGenerationEU (PRIN 2022 project ASTRA, Grant no. 20228YYR5Z\_002,
			CUP J53D23000640006).}
}

\author{Xiaoxing Ren   \and
Michelangelo Bin \and
        Ivano Notarnicola 		\and
        Thomas~Parisini 
}

\institute{
X. Ren \at
Department of Electrical and Electronic Engineering, Imperial College London, London, UK\\ 
            \email{xiaoxing.ren@imperial.ac.uk}
            \and
            M. Bin \at
           	Department of Electrical, Electronic, and Information Engineering, University of Bologna, Italy\\ 
            \email{michelangelo.bin@unibo.it}           
           \and
           I. Notarnicola\at
           	Department of Electrical, Electronic, and Information Engineering, University of Bologna, Italy\\ 
           	\email{ivano.notarnicola@unibo.it}   
           \and 
           T. Parisini \at 
           Department of Electrical and Electronic Engineering, Imperial College London, London, UK
           	\at Department of Electronic Systems, Aalborg University, Denmark
           	\at 
           	Department of Engineering and Architecture, University of Trieste, Trieste, Italy\\ 
           	\email{t.parisini@imperial.ac.uk}
}
\date{Received: date / Accepted: date}

\maketitle

\begin{abstract}
This paper studies a distributed algorithm for  constrained consensus optimization that is obtained by fusing the Arrow-Hurwicz-Uzawa primal-dual gradient method for centralized constrained optimization and the Wang-Elia method for distributed unconstrained optimization.
It is shown that the optimal primal-dual point is a semiglobally exponentially stable equilibrium for the algorithm, which implies linear convergence.
The analysis is based on the  separation between a slow centralized optimization dynamics describing the evolution of the average estimate toward the optimum, and a fast   dynamics describing the evolution of the consensus error over the network. These two dynamics are mutually coupled, and the stability analysis builds on control theoretic tools such as time-scale separation, Lyapunov theory, and the small-gain principle.
Our analysis approach highlights that the consensus dynamics can be seen as a fast, parasite one, and that stability of the distributed algorithm is obtained as a robustness consequence of the semiglobal exponential stability properties of the centralized method. This perspective can be used to enable other significant extensions, such as time-varying networks or delayed communication, that can be seen as ``perturbations" of the centralized algorithm.

\keywords{Convex programming \and Distributed Constrained Optimization \and Discrete-Time Optimization Algorithms \and Stability Analysis of Optimization Algorithms.}
\subclass{65K05 \and 93A14 \and 93A16 \and 90C25 \and 90C30 \and 90C35 \and 93D05 \and 93D23}
\end{abstract}


\section{Introduction}

\subsection{Problem overview} \label{problem_formulation}
Given a set $\cN\coloneqq \{1,\dots,n\}$ of $n\in\N$ agents  connected through a network, we consider the problem of designing a distributed algorithm by which agents cooperatively solve a constrained consensus optimization problem of the following form
\begin{equation}\label{d.optimization_problem}
	\begin{aligned}
		\min_{\theta\in\R} \: & \: \sum_{i\in\cN} f_i(\theta) \, ,
		\\
		\subj \: & \: g_i(\theta)\le 0,\ i=1,\dots,n \, ,
	\end{aligned}
\end{equation} 
in which ``$\le$" is applied component-wise and, for each $i\in\cN$, $f_i:\R\to\R$ and $g_i:\R\to\R^{m_i}$ ($m_i\in\N$) are ``local''  functions only accessible to agent $i$. 
In~\eqref{d.optimization_problem} we assume that $\theta$ is a scalar (i.e., $\theta\in\R$) to ease the notation and the technical derivations; however, the presented arguments easily extend to cover the case where $\theta\in\R^{n_\theta}$ for some arbitrary $n_\theta\in\N$. 
Under the assumption that each agent $i$ can evaluate the gradient of $f_i$ and $g_i$ at an arbitrary point and can communicate with a subset $\cN_i\subset\cN\setminus\{i\}$ of neighboring agents, this paper studies the following {\it discrete-time} algorithm 
\begin{equation}\label{s.algorithm_componentwise}
	\begin{aligned}
x_i^{t+1} &= x_i^t + \sum_{j\in\cN_i} k_{ij} (x_j^t-x_i^t + z_j^t -z_i^t) - \gamma{n} 
\nabla \ell_i(x_i^t, \lambda_i^t) \, , \\
z^{t+1}_i &= z_j^t - \sum_{j\in\cN_i} k_{ij} (x_j^t-x_i^t)\, ,\\
\lambda_i^{t+1} & = \max\left\{0,\lambda_i^t+\gamma g_i(x_i^t)  \right \}, \qquad \lambda_i^0\ge 0\, ,
	\end{aligned}
\end{equation}
where $x_i\in\R$ denotes the estimate agent $i$ has on the optimal solution $\theta\sr$ of Problem~\eqref{d.optimization_problem}, $z_i\in\R$ and $\lambda_i\in(\mathbb{R}_{\geq 0})^{m_i}$ are auxiliary state variables, $\nabla\ell_i$ denotes the gradient of
\begin{equation*}
	\ell_i(x_i^t,\lambda_i^t) \coloneqq  f_i(x_i^t) +  \sprod{\lambda_i^t}{g_i(x_i^t)}  
\end{equation*}
with respect to the first variable,  and  $\gamma>0,\ k_{ij} >0$ are suitable design parameters.  
Algorithm~\eqref{s.algorithm_componentwise} can be interpreted as ``mixing" the distributed unconstrained gradient flow proposed in~\cite{wang2011control} and studied in~\cite{binStabilityLinearConvergence2022} with the centralized Arrow-Hurwicz-Uzawa  primal-dual algorithm for constrained problems dealt with in~\cite{bin2024semiglobal}. 


%

While discrete-time algorithms are better suited for practical applications, their convergence analysis presents significant challenges. 
In turn, to the best of the authors' knowledge, a proof of exponential convergence for~\eqref{s.algorithm_componentwise} and affine discrete-time algorithms remains an open problem, even  if $\sum_{i\in\cN}f_i$ is strongly convex  and $g$ is convex (which shall be assumed later on). This article aims to address this gap by providing a purely discrete-time stability analysis for~\eqref{s.algorithm_componentwise}. Specifically, it is shown that the primal-dual solution of \eqref{d.optimization_problem} is a Lyapunov stable and semiglobally exponentially attractive equilibrium of Algorithm~\eqref{s.algorithm_componentwise}. As a simple counterexample shows (see~\cite{bin2024semiglobal}), this is the best result achievable for such an algorithm since global exponential stability cannot be achieved in general.

\subsection{Related Works}



Distributed optimization problems with local constraints, such as~\eqref{d.optimization_problem}, have been extensively studied in the literature. Most works focus on continuous-time algorithms, owing to the mature development of continuous-time stability theory and the fact that the stability analysis is easier in continuous time~\cite{wang_control_2010,liuCollectiveNeurodynamicApproach2017,zhu_projected_primal_dual,tianDistributedOptimizationMultiagent2021,jiangContinuousTimeAlgorithmApproximate2021,liDistributedContinuousTimeNonsmooth2020,Ren_Distributed_Global_Optimization,cisneros2021distributed,li2020smooth,liExponentiallyConvergentAlgorithm2022,zhu_projected_primal_dual}.
For these algorithms, linear convergence is not typically guaranteed.
An exception is the algorithm  proposed in \cite{zhu_projected_primal_dual}, which achieves linear convergence although only under linear equality constraints.
Local linear equality and inequality constraints are addressed in~\cite{wang_control_2010}, which also provides a control theoretic perspective on the problem. These works have been then extended by~\cite{hatanaka_passivity-based_2018,yamashita2020passivity} to address the case of local convex inequality constraints; the stability analysis is grounded on passivity theory.
Continuous algorithms generally lend themselves to a simpler analysis and, consequently, to more powerful results. However, they need to be discretized to be used, and their discretization presents significant challenges often hindering their applicability.
In particular, most continuous-time algorithms solving constrained problems present discontinuities, and continuity is instead key to apply basic discretization theorems (see, e.g., \cite[Sec.~2]{stetter_analysis_1973}). Moreover, discretization methods leading to first-order algorithms, favorable in applications, use discretization steps that in general depend on the integration interval, which can be taken uniform over an infinite time horizon only under very specific conditions (e.g., exponential contraction, or asymptotic stability of a compact set~\cite[Thm. 2]{varagnolo_newton-raphson_2016}) and need  to diminish  in time otherwise. Such conditions, however, are  typically not proved, and when they are, other phenomena may occur such as loss of globality of convergence \cite{bin2024semiglobal}. We also highlight that diminishing stepsizes lead to sublinear convergence and loss of robustness, the latter that is key to handle perturbation phenomena such as time-varying communication networks and delays.
%
%

Purely discrete-time algorithms are therefore generally preferable. However, they require a more involved analysis. 
%
%
A substantial number of existing works exist in the literature.  Reference~\cite{nedic2010constrained} proposes a purely primal approach to constrained optimization, which converges under a diminishing stepsize without any convergence rate guarantees. A diminishing stepsize is employed also in~\cite{changDistributedConstrainedOptimization2014}, where a consensus-based discrete-time primal-dual perturbation method is proposed to solve convex optimization problems with coupled inequality constraints, and in~\cite{liuDistributedDiscreteTimeAlgorithms2021}, which investigates convex optimization problems with local convex inequality constraints and linear equality constraints under weight-unbalanced topology. The convergence rate is sublinear due to the diminishing stepsize.
Time-varying unbalanced networks are considered in
\cite{xieDistributedConvexOptimization2018}, where a discrete-time algorithm is proved to asymptotically converge to the optimal solution of a convex constrained problem. Moreover,
\cite{notarnicolaConstraintCoupledDistributedOptimization2020} and \cite{camisaDistributedStochasticDual2022} consider distributed optimization problems with local set constraints and coupling inequality constraints. They propose discrete-time algorithms in which each agent iteratively performs a constrained, small-sized local optimization followed by a distributed dual update.
An augmented Lagrangian tracking-based algorithm is proposed in~\cite{falsoneAugmentedLagrangianTracking2023} for distributed optimization problems coupled by linear equality and nonlinear inequality constraints. The algorithms needs solving an intermediate local minimization problem at each iteration.
A discrete-time primal-dual method that can deal with time-varying directed networks is proposed in
\cite{gongPrimalDualAlgorithmDistributed2024} for convex optimization problems with coupled affine equality and (not necessarily affine) inequality constraints over time-varying and directed networks.  
Reference~\cite{mateos2016distributed} proposes instead a distributed saddle-point subgradient algorithm  with Laplacian averaging.
Finally,~\cite{liDistributedProximalAlgorithms2021} proposes a proximal primal-dual algorithm with convergence rate $\mathcal{O}(1/\sqrt{t})$ to deal with convex problems with local set constraints and coupled inequality constraints,  while~\cite{wuDistributedOptimizationCoupling2023} proposes an integrated primal-dual proximal algorithm with convergence rate $\mathcal{O}(1/t)$ for solving convex problems with both coupled inequality and linear equality constraints.

Although linear convergence has been extensively studied for discrete-time unconstrained algorithms~\cite{shiLinearConvergenceADMM2014,shiEXTRAExactFirstOrder2015,nedicAchievingGeometricConvergence2017,quHarnessingSmoothnessAccelerate2018,quAcceleratedDistributedNesterov2020,puPushPullGradient2021,alghunaimDecentralizedProximalGradient2021,binStabilityLinearConvergence2022,renAcceleratedDistributedGradient2022,notarnicola_gradient_2023}, and results also exist for problems with linear equality constraints (see, e.g.,~\cite{andersonDistributedApproximateNewton2019}), 
to the best of the author's knowledge, proving stability and exponential convergence for discrete-time algorithms handling inequality constraints is still an open problem, even under the assumptions that $\sum_{i\in\cN} f_i$ is strongly convex and $g_i$ is convex for each $i\in\cN$.
%



\subsection{Contribution} 
This paper proves semiglobal exponential stability of the optimal primal-dual solution of the constrained problem \eqref{d.optimization_problem} for Algorithm~\eqref{s.algorithm_componentwise}.
As commented in the previous section, this is an open problem, not only for what concerns Algorithm~\eqref{s.algorithm_componentwise}, but for discrete-time first-order methods in general. In the specific, under strong convexity of the cost function $\sum_{i\in\cN}f_i$ and convexity of the constraint functions $g_i$, it is proved that the unique optimal primal-dual solution of Problem~\eqref{d.optimization_problem} is a Lyapunov stable equilibrium of~\eqref{s.algorithm_componentwise} and, for every compact set of initial conditions for~\eqref{s.algorithm_componentwise}, there exists a sufficiently small upper-bound for the stepsize $\gamma$, such that every sequence generated by~\eqref{s.algorithm_componentwise} originating in the initialization set exponentially converges to such an optimal equilibrium.
We stress that this is the best convergence result achievable for this class of algorithms, as the counterexample provided in~\cite{bin2024semiglobal} for the centralized version shows.

The stability analysis builds on control tools such as time-scale separation, Lyapunov theory, and the small-gain principle. In particular, the dynamics of \eqref{s.algorithm_componentwise} is split into a  slow centralized dynamics solving a centralized constrained optimization problem, and a fast network dynamics describing the evolution toward consensus.
The centralized dynamics coincides with the Arrow-Hurwicz-Uzawa primal-dual gradient method analyzed in \cite{bin2024semiglobal}, which is proved there to be semiglobally exponentially stable. The robustness properties of exponential stability are exploited here to deal with the fast network dynamics, which is therefore treated as a perturbation of the centralized algorithm that can be overcome by a suitably small stepsize.
While specialized to handle such a network consensus phenomena, the underlying rationale in which perturbations are handled thanks to the robustness margin of exponential stability can be easily extended to cover other problems of interest, such as time-varying or time-delayed communications, whenever such phenomena can be reduced to fast parasite dynamics that become negligible when the stepsize is sufficiently small. This analysis approach is enabled by the fact that Algorithm~\eqref{s.algorithm_componentwise} employs a non-diminishing stepsize.

As a byproduct of the stability analysis, an explicit theoretical upper bound for the stepsize to guarantee semiglobal linear convergence is provided. Such a bound highlights how the network size, structure, and connectivity, the Lipschitz constant, and the convexity properties of the involved functions  influence the stepsize and the convergence rate.  
Finally, we remark that, for simplicity of exposition,  a single stepsize for both the primal and the dual updates is assumed. While this is generally a more difficult case to analyze, as it prevents a further time-scale separation between the primal and the dual dynamics (see, e.g.,~\cite{alghunaim2020linear}), there are cases where different stepsize values are preferable (for instance, if one wants to embed the multiplication by $n$ of $\nabla\ell_i(x_i^t,\lambda_i^t)$ in \eqref{s.algorithm_componentwise} within the stepsize bound). The analysis carried out here can be easily extended to cover such a case.

\subsection{Notation} \label{notation} 
If $\sim$ is a binary relation on a set $S$  and $z\in S$, we let $S_{\sim z}\coloneqq\{ s\in S\st s\sim z\}$.
If not otherwise specified, binary relations are applied component-wise to vectors. 
We denote by $\nabla f$ the gradient of a differentiable function $f$;
 $\sigma(A)$ denotes the spectrum of a matrix $A\in \R^{n\times n}$; $A$ is  Schur if $\sigma(A)$ is contained in the open unit disk in the complex plane; $A>0$ means that $A$ is positive definite.
If $\cI$ is an ordered set of finite cardinality $n$,  $(x_i)_{i\in\cI}$ denotes the $n$-tuple of $\R^n$ ordered by $\cI$ in the obvious way. 
Set inclusion (either strict or not) is denoted by $\subset$.
With $n\in\N$, we let $\1\coloneqq (1,\dots,1)\in\R^n$, $I \in \R^{n \times n}$ be the identity matrix, and
$\sprod{x}{y}\coloneqq \sum_{h=1}^n x_h y_h$ denote the standard inner product of two vectors $x, y \in \R^n$. Moreover,
$| \cdot |$ denotes the Euclidean norm of a vector or the matrix-induced 2-norm, and the closed ball of radius $r$ around point $\bar{x}\in \R^n$ is denoted by 
$\overline{\mathbb{B}}_r(\bar x)\coloneqq \{ x \in \R^n \st |x-\bar{x}|\le r \}$.

We let $S\in\R^{n\x(n-1)}$ be a matrix satisfying 
\begin{align}\label{d.S}
	S\T \1 &=0, & S\T S &= I_{n-1}.
 \end{align}
Moreover, we define the matrix $T\in\R^{n\x n}$ and its inverse as
\begin{equation}\label{d.T}
	T: =\begin{bmatrix}
		\frac{\1\T}{n} \\[.5ex] S\T 
	\end{bmatrix}, \qquad T\inv : = \begin{bmatrix}
		\1& S
	\end{bmatrix}.
\end{equation}
From \eqref{d.S}-\eqref{d.T}, we deduce that  
\begin{align}\label{e.I_one_S}
	I_n &= \tfrac{1}{n}\1\1\T + SS\T, & |S|&=1.
\end{align}
Given $\chi\in\R^n$, its \emph{average-dispersion decomposition}~\cite{binStabilityLinearConvergence2022} is the pair $(\chi\m,\chi\p)=T\chi$, where $\chi\m\coloneqq \1\T\chi/n\in\R$  and $\chi\p\coloneqq S\T \chi\in\R^{n-1}$ are called, respectively, the \emph{average} and the \emph{dispersion} components of $\chi$. 
From \eqref{d.S} and~\eqref{d.T}, we obtain  $\chi = T\inv(\chi\m,\chi\p)= \1 \chi\m + S\chi\p$, and $|\chi|^2 = n\chi\m^2+|\chi\p|^2$.

\section{Main Result}
This section presents the main result of the paper. First, Section~\ref{sec.main.ass} details the standing assumptions. Subsequently,   Section~\ref{sec.main.equilibria}  links the optimal solution of Problem~\eqref{d.optimization_problem} to the equilibria of Algorithm~\eqref{s.algorithm_componentwise}. Finally,   Section~\ref{sec.main.result} states the main stability result.

\subsection{Standing Assumptions} \label{sec.main.ass}
Algorithm \eqref{s.algorithm_componentwise} is studied under a number of assumptions reported hereafter.
First, rather standard convexity, smoothness, and feasibility properties are assumed for Problem \eqref{d.optimization_problem}.
\begin{assumption}\label{ass.optimization_problem_data} 
The function {$\sum_{i\in\cN}f_i$} is strongly convex and twice continuously differentiable. Moreover, for each $i\in\cN$, function $g_i$ is convex and twice continuously differentiable. Finally, there exists $\bar\theta\in\R$ such that, for every $i\in\cN$, $g_i(\bar\theta)\le 0$.
\end{assumption}

Under Assumption~\ref{ass.optimization_problem_data}, Problem~\eqref{d.optimization_problem} admits a unique optimal solution $\theta\sr\in\R$ \cite[Lem.~1]{bin2024semiglobal}.
Let
\begin{equation*}
	A(\theta\sr) \coloneqq \{ (i,j)\in\cN\x\N\st j\in\{1,\dots,m_i\} ,\    g_{i,j}(\theta\sr) = 0 \}
\end{equation*}
denote the set of all indices corresponding to the constraints that are  \emph{active} at the optimum $\theta\sr$.
Then, in addition to Assumption~\ref{ass.optimization_problem_data}, we assume that $\theta\sr$ is \emph{regular} in the following sense \cite[Sec. 3.1]{bertsekas1999nonlinear}.

\begin{assumption}\label{ass.regular_θ*}
	The vectors $\{ \nabla g_{i,j}(\theta\sr)\st (i,j)\in A(\theta\sr)\}$ are linearly independent.
\end{assumption}

Assumptions~\ref{ass.optimization_problem_data} and~\ref{ass.regular_θ*} imply the existence, for each $i\in\cN$, of a \emph{unique}  $\lambda\sr_i\in\R^{m_i}$ such that the following \emph{Karush-Kuhn-Tucker (KKT) conditions} hold \cite[Prop.~3.3.1]{bertsekas1999nonlinear}
\begin{subequations}\label{d.KKT}
	\begin{align}
		&\sum_{i\in\cN}  \nabla \ell_i(\theta\sr,\lambda\sr_i) = 0 \, , \label{d.KKT.gradients}\\
		&\lambda\sr_i \ge 0,\ g_i(\theta\sr)\le 0, && \forall i\in\cN \, ,  \label{d.KKT.constraints}\\
		&\lambda\sr_{i,j} g_{i,j}(\theta\sr) = 0,&&\forall i\in\cN \, , \ \forall j=1,\dots,m_i  \, .\label{d.KKT.active_constraints}
	\end{align}
\end{subequations}
Conditions~\eqref{d.KKT} are also sufficient: If $(\theta,\lambda)$ satisfies \eqref{d.KKT}, then $\theta=\theta\sr$ is the optimal solution of \eqref{d.optimization_problem} \cite[Prop.~3.3.4]{bertsekas1999nonlinear}.

The topology of the network through which agents exchange information  is completely determined by the \emph{neighborhoods} $\cN_i\subset \cN \setminus\{i\}$  introduced in Section \ref{problem_formulation}.  Regarding the communication weights, the following assumption is introduced.
\begin{assumption}\label{ass.network}  
The coefficients $k_{ij}>0$ of \eqref{s.algorithm_componentwise} are such that the matrix
$K\in\R^{n\x n}$, whose $(i,j)$th element $K_{ij}$   is defined as 
\begin{equation*}
	K_{ij} \coloneqq 
	\begin{cases}
		\sum_{j\in\cN_i} k_{ij} \, , & \text{if}\ j=i \, ,\\
		-k_{ij} \, , & \text{if}\ j\in\cN_i \, , \\
		0\, , & \text{otherwise},
	\end{cases}
\end{equation*} 
satisfies the following conditions
\begin{align}\label{d.K}
	K&=K\T, & \ker K &= \linspan\1, & \sigma(K)&\subset{\left[0,  1 \right).}
\end{align}
\end{assumption}
Assumption~\ref{ass.network} (which is the same made in \cite{binStabilityLinearConvergence2022} for the unconstrained case) implies that the communication network is undirected and strongly connected.
Moreover,~\eqref{d.K} implies that $\1\T K=0$ and $S\T KS$ invertible.

\subsection{Optimality and Equilibria}\label{sec.main.equilibria}
Algorithm \eqref{s.algorithm_componentwise} can be rewritten in compact form as
\begin{subequations}\label{s.algorithm}
	\begin{align}
x^{t+1} &= (I-K) x^t - K z^t - \gamma {n}\nablaL(x^t, \lambda^t) \, ,
\label{s.algorithm.x}\\
z^{t+1} &= z^t + K x^t \, ,
\label{s.algorithm.z}\\
		\lambda^{t+1} &= \max\{0,\,\lambda^t+\gamma g(x^t)\} \,, \qquad \lambda^0 \in(\R_{\ge 0})^m \, ,
		\label{s.algorithm.λ}
	\end{align}
\end{subequations} 
in which $K$ is the matrix defined in Assumption~\ref{ass.network},  $x\coloneqq (x_i)_{i\in\cN} \in\R^n$, $z\coloneqq (z_i)_{i\in\cN} \in\R^n$,  
 $\lambda\coloneqq (\lambda_i)_{i\in\cN} \in(\R_{\ge 0})^m ${, with $m\coloneqq m_1+\cdots+m_n$}, $g :\R^n\to\R^m$ is defined as $g(x)\coloneqq (g_i(x_i))_{i\in\cN}$, and $\nablaL  :\R^n\x(\R_{\ge 0})^m \to \R^n$ is defined as $\nablaL(x,\lambda) \coloneqq  ( \nabla\ell_i(x_i,\lambda_i) )_{i\in\cN}$.

The following lemma characterizes the equilibria of \eqref{s.algorithm} in terms of the optimality conditions \eqref{d.KKT}, and it is a generalization of \cite[Lem.~2]{bin2024semiglobal}, which applies to the constrained but centralized case, and \cite[Lem.~1]{binStabilityLinearConvergence2022}, which applies to the unconstrained but distributed case.
We underline that the result of the lemma uses only part of Assumptions
\ref{ass.optimization_problem_data} and \ref{ass.regular_θ*}, i.e., differentiability of $\sum_{i\in\cN}f_i$ and $g_i$. 
\begin{lemma} \label{lem.equilibria}
	Suppose that Assumption~\ref{ass.network} holds and that $\sum_{i\in\cN}f_i$ and $g_1,\dots, g_n$ are differentiable. Then, every equilibrium $(x\eq,z\eq,\lambda\eq)$ of \eqref{s.algorithm} satisfies $x\eq = \1 \theta\sr$ and $\lambda\eq=\lambda\sr$, where $(\theta\sr,\lambda\sr)$ is a solution of~\eqref{d.KKT}. 
	Conversely, if $(\theta\sr,\lambda\sr)$ satisfies \eqref{d.KKT}, there exists an equilibrium $(x\eq,z\eq,\lambda\eq)$ of \eqref{s.algorithm} such that  $x\eq=\1\theta\sr$ and $\lambda\eq=\lambda\sr$.
\end{lemma}
\begin{proof}
	See Appendix~\ref{proof.lem.equilibria}.
\end{proof}
Lemma \ref{lem.equilibria}  implies that system \eqref{s.algorithm} has an entire closed subspace $\cA\sr$ of equilibria, all of which are optimal points for \eqref{d.optimization_problem} in the sense that they satisfy \eqref{d.KKT}. In particular, as clear from the proof of the lemma, the set $\cA\sr$ is given by
\begin{align*}
\cA\sr \coloneqq  \big\{ & (x\eq,z\eq,\lambda\eq)\in\R^{n}\x\R^{n}\x(\R_{\ge 0})^{m}\st \\
&\qquad \exists (\theta\sr,\lambda\sr) \ \text{solving \eqref{d.KKT}}, \ x\eq = \1 \theta\sr,\ \lambda\eq=\lambda\sr,   \\
&\qquad z\eq \in   -\gamma{n} S (S\T K S)\inv S\T \Phi(\1\theta\sr,\lambda\sr) + \linspan\1 \big\} \, .
\end{align*}
The set $\cA\sr$ is closed, but not compact.  
Moreover, if a solution of \eqref{s.algorithm} converges to $\cA\sr$, the actual equilibrium point to which it converges depends on the initial conditions.
Indeed, it is an immediate consequence of \eqref{d.K} that the quantity $\1\T z^t$ is constant along every solution of \eqref{s.algorithm}. Since we can write $z^t = (SS\T + \frac{\1\1\T}{n}) z^t$, then it can be concluded that every solution of \eqref{s.algorithm} that converges to $\cA\sr$,   converges to the point $(\1\theta\sr, z\eq, \lambda\sr)$, where $(\theta\sr,\lambda\sr)$ satisfies \eqref{d.KKT} and $z\eq=  -\gamma {n}S (S\T K S)\inv S\T \Phi(\1\theta\sr,\lambda\sr) + \1 \frac{1}{n} \sum_{i\in\cN} z_i^0$. In view of this, our focus will not be the distance of $(x,z,\lambda)$ to one particular equilibria of $\cA\sr$, but the distance of $(x,z,\lambda)$ to the whole set $\cA\sr$, which is defined as $|(x,z,\lambda)|_{\cA\sr} \coloneqq  \inf_{a\in\cA\sr} |(x,z,\lambda) - a|$.

\subsection{Stability Result}\label{sec.main.result}
The theorem reported hereafter establishes the main result of the paper, namely, semiglobal exponential stability of $\cA\sr$ for Algorithm~\eqref{s.algorithm_componentwise}.
Specifically, it is shown that, for every arbitrarily large compact set of initial conditions, there exists a maximum value for the parameter $\gamma$ such that, for all  positive values of $\gamma$ smaller than such a maximum value, the set $\cA\sr$ of optimal equilibria is exponentially asymptotically stable for Algorithm~\eqref{s.algorithm_componentwise} with a domain of attraction that includes the chosen set of initial conditions. 

\begin{theorem}\label{thm.main}
Suppose that Assumptions~\ref{ass.optimization_problem_data}, \ref{ass.regular_θ*}, and \ref{ass.network} hold.
Let $\theta\sr\in\R$ be the unique solution of \eqref{d.optimization_problem} and $\lambda\sr\in\R^m$ the unique {vector} such that $(\theta\sr,\lambda\sr)$ solves \eqref{d.KKT}. 
Then, for every compact subset $\Xi_0\subset \R^n\x\R^n\x(\R_{\ge 0})^m$ of initial conditions for \eqref{s.algorithm}, there exists $\bar\gamma>0$, and for every $\gamma\in(0,\bar\gamma)$, there exist $\mu=\mu(\gamma)\in(0,1)$ and $c = c(\gamma)\ge 0$, such that every solution $(x,z,\lambda)$ of Algorithm~\eqref{s.algorithm}, with $(x^0,z^0,\lambda^0)\in\Xi_0$, satisfies
	\begin{equation}\label{e.exp-bound}
		\forall t\in\N,\quad |(x^t,z^t,\lambda^t)|_{\cA\sr} \le c \mu^t |(x^0,z^0,\lambda^0)|_{\cA\sr} \, .
	\end{equation}
\end{theorem}

The proof of Theorem~\ref{thm.main} is provided in Section~\ref{sec.proof.main-theorem}.
Theorem~\ref{thm.main} implies several properties for Algorithm~\eqref{s.algorithm}. First, provided that $\gamma$ is sufficiently small, all trajectories generated by the algorithm from a given arbitrary initialization set $\Xi_0$ are equibounded, i.e., the union of the orbits of all solutions of \eqref{s.algorithm} from $\Xi_0$ is a precompact set. This property is also known as \emph{Lagrange stability on $\Xi_0$} \cite{andriano_global_1997}.
Second, since $\mu\in(0,1)$, 
Theorem \ref{thm.main} also implies \emph{Lyapunov stability of $\cA\sr$ relative to $\Xi_0$}, i.e., for every neighborhood $V$ of $\cA\sr$, there exists  a neighborhood $U$ of $\cA\sr$, such that every solution of \eqref{s.algorithm} originating inside $U\cap \Xi_0$ satisfies $(x^t,z^t,\lambda^t)\in V$, for all $t\ge 0$. This is a continuity property that implies that ``small'' deviations of the initial conditions from the target steady-state set will lead to ``small'' fluctuations around it.
Finally, Theorem~\ref{thm.main} implies convergence of all the agent sequences $x_i^t$ generated by Algorithm~\ref{s.algorithm} from $\Xi_0$ to the optimal point $\theta\sr$ (the solution to the optimization problem~\eqref{d.optimization_problem}) and of $\lambda^t$ to the Lagrange multiplier $\lambda\sr$ such that $(\theta\sr,\lambda\sr)$ solves the KKT conditions~\eqref{d.KKT}.
Convergence, as that of $(x^t,z^t,\lambda^t)$ to $\cA\sr$, is exponential with a rate $\mu$ that depends on $\gamma$.
Specifically, as clear from the proof (see, in particular, \eqref{d.omega_gamma}, \eqref{d.T-eta}, \eqref{d.proof.eta2}, and \eqref{d.proof.c_mu}), the constants $c$ and $\mu$ are estimated as
\begin{align*}
	\mu &=  \sqrt{1-\omega} \, ,  & c &= n^{{\frac{1}{2}}} \sqrt{\max\left\{3, \frac{c_u}{c_l}\right\}} \left(\frac{1}{1-\omega}\right)^{\frac{T}{2}},
\end{align*}
in which
\begin{align*}
	\omega & \coloneqq  \min\left\{ \frac{1}{12}\gamma \mu_f, \, \frac{1}{12}\gamma^2 k_2^2,\, \frac{1}{6}\gamma\frac{h}{\kappa_1},\, \frac{1}{6},\, \frac{1}{8c_u}  \right\} \in (0,1) \, ,\\
	T &\coloneqq \frac{3\kappa_0^2 - \min\left\{1,\gamma\frac{\mu_f}{4\omega} \right\}\varepsilon^2 }{ \omega \min\left\{1,\gamma\frac{\mu_f}{4\omega} \right\}\varepsilon^2}\, .
\end{align*}
In the previous definitions, $c_u$ and $c_l$ (defined in Section~\ref{sec.proof.preliminary}) are the maximum and minimum eigenvalues of the solution $P$ of the Lyapunov equation characterizing the network's consensus dynamics. In particular the ratio $c_u / c_l$ coincides with the condition number of $P$ and it is related to the connectivity, hence stability, of the network;  less connected and stable networks are associated with a larger ratio  $c_u / c_l$, hence with a larger value of $c$.
The expression of $c$ also reveals that the exponential bound of Theorem~\ref{thm.main} is not independent from the network size, as $c$ grows with $n$. It is currently unknown if such a bound can be made independent from $n$.
The constant $\mu_f$ appearing in the definition of $\omega$ and $T$ is the convexity parameter of the cost function  $\sum_{i\in\cN}f_i$ (see \eqref{d.mu_f}); the constant $\kappa_0$ is the radius of a closed ball containing $(x\m,\lambda,x\p,z\p)$ for every $(x,z,\lambda)\in \Xi_0$ (see \eqref{d.kappa0-eq}); $\kappa_1>\kappa_0$ is the radius of a larger ball defined in \eqref{d.kappa1}; $k_2$ is the Lipschitz constant of the Lagrangian $\sum_{i\in\cN} \ell_i$ on a compact superset $\cK$ of the previously-defined balls formally introduced in~\eqref{set_k}; $h>0$ is defined as the minimum absolute value of the inactive constraint functions at the optimal point; finally, $\varepsilon>0$ is a small number formally defined by \eqref{d.epsilon}.
Notice that, for small enough $\gamma$, one can assume $\omega=\frac{1}{12}\gamma^2 k_2^2<1$ (the last inequality follows from \eqref{gamma_all} and by the definition of $\bar\gamma_{16}$ in \eqref{d.bargamma_2}). Hence, $\mu \ge \sqrt{11/12}$,  which gives a uniform (with respect to $\gamma$) lower bound on the convergence rate for small values of~$\gamma$.

\begin{remark}
	It is worth remarking that a \emph{global} convergence result, in which the maximum value of $\gamma$ is the same for every initialization set, cannot be achieved even with sub-linear convergence and with an assumption of Lipschitz gradients, as the simple counterexample in~\cite{bin2024semiglobal} shows. Hence, semiglobal convergence is the best property one can achieve in the considered case. 
	Moreover, one can attribute this loss of globality to the presence of constraints, and not to decentralization. Indeed, global exponential stability holds for the unconstrained version of the algorithm under the assumption of Lipschitz gradients~\cite{binStabilityLinearConvergence2022}, while it fails for the centralized version with constraints. 
	As the analysis idea behind the proof of the theorem is that the average component of Algorithm~\eqref{s.algorithm} behaves like the centralized algorithm, perturbed by the network effect due to the consensus dynamics, the properties of the centralized algorithm are necessarily inherited by its distributed version. 
\end{remark} 

\section{Stability Analysis}\label{sec.proof.main-theorem}
%
This section provides the proof of Theorem~\ref{thm.main}.
The proof follows a canonical nonlinear control paradigm based on time-scale separation, Lyapunov analysis, and the small-gain principle. In the proof, we highlight the time-scale separation properties underlying the functioning of the algorithm~\eqref{s.algorithm}, and the role of control theory in establishing stability and convergence. 
Throughout the proof, we shall suppose that all the assumptions of the theorem hold, and we fix now, once and for all, an arbitrary compact set $\Xi_0\subset\R^n\x\R^n\x(\R_{\ge 0})^m$ of initial conditions for~\eqref{s.algorithm}.

\subsection{The Core Subsystem}\label{sec.proof.core_subsystem}
This section follows the argument of \cite{binStabilityLinearConvergence2022} to construct a  ``core subsystem'' obtained by eliminating from \eqref{s.algorithm} the marginally stable dynamics of \eqref{s.algorithm.z}. The resulting system has an optimal equilibrium whose asymptotic stability implies asymptotic stability of $\cA\sr$ for the original algorithm \eqref{s.algorithm}. 

Consider the average-dispersion decomposition (Section \ref{notation}) of the variables $x$ and $z$, namely 
\begin{align*}
	(x\m,x\p) &\coloneqq  Tx, &(z\m,z\p) &\coloneqq  Tz.
\end{align*}
In view of \eqref{d.K},   every solution of \eqref{s.algorithm} satisfies
\begin{equation}\label{s.zm}
	\forall t\in\N,\quad z\m^t = z\m^0.
\end{equation}
Let ($\theta\sr,\lambda\sr)$ be the unique solution of \eqref{d.KKT} such that $\theta\sr$ solves \eqref{d.optimization_problem}, and let 
\begin{equation}\label{d.z_sr}
	 z\p\sr\coloneqq -\gamma {n} (S\T K S)\inv S\T\Phi(\1\theta\sr,\lambda\sr).
\end{equation} 
Then,  
\begin{equation}\label{e.dist_A_without_zm}
\begin{aligned}
	|(x,z,\lambda)|_{\cA\sr} &= \inf_{(x\eq,z\eq,\lambda\eq)\in\cA\sr} |(x-x\eq,z-z\eq,\lambda-\lambda\eq)|
	\\
	&= \inf_{c\in\R} |(x-\1\theta\sr,z-(Sz\p\sr + \1 c),\lambda-\lambda\sr)|\\
	&= | (\sqrt{n}(x\m-\theta\sr), x\p, z\p-z\p\sr, \lambda-\lambda\sr)|,
\end{aligned} 
\end{equation}
where the last equality follows from the fact that
\begin{align*}
  & |(x-\1\theta\sr,z-(Sz\p\sr + \1 c),\lambda-\lambda\sr)| 
\\
&\ = |(\1(x\m-\theta\sr)+Sx\p, S(z\p-z\p\sr)+\1(z\m-c),\lambda-\lambda\sr)|	\\&\ = \Big( |\1(x\m-\theta\sr)+Sx\p|^2 + |S(z\p-z\p\sr)+\1(z\m-c)|^2  + |\lambda-\lambda\sr|^2 \Big)^{\frac{1}{2}}
\\
&\ = \Big( n |x\m-\theta\sr|^2 + |x\p|^2 + n |z\m-c|^2 + |z\p-z\p\sr|^2 + |\lambda-\lambda\sr|^2 \Big)^{\frac{1}{2}} \\
&\ = | (\sqrt{n}(x\m-\theta\sr), x\p, \sqrt{n}(z\m-c), z\p-z\p\sr, \lambda-\lambda\sr )| \, ,
\end{align*}
which has a minimum for $c=z\m$.

Furthermore, in view of \eqref{d.K}, $Kz = K(Sz\p + \1 z\m) = KSz\p$, so as the variable $z\m$ does not affect   any other state variable.
Since \eqref{e.dist_A_without_zm} implies that $z\m$ has no effect on $|(x,z,\lambda)|_{\cA\sr}$ as well, and since  \eqref{s.zm}    implies that $z\m$ is marginally stable,  we  can therefore neglect the dynamics of $z\m$ and only focus on that of $(x,z\p,\lambda)$, which, in the average-dispersion coordinates, is described by  
\begin{subequations}\label{s.reduced}
	\begin{align}
		x\m^{t+1} &= x\m^t - \gamma\1\T \nablaL(x^t,\lambda^t)
		\label{s.reduced.xm} \, , \\
		x\p^{t+1} &= (I-S\T KS) x\p^t - S\T KS z\p^t  - \gamma{n} S\T    \nablaL(x^t,\lambda^t)  \, , 
		\label{s.reduced.xp}\\
		z\p^{t+1} &= z\p^t + S\T K S x\p^t \, , 
		\label{s.reduced.zp} \\
		\lambda^{t+1} &= \max\{0,\, \lambda^t + \gamma g(x^t)\} \, .
		\label{s.reduced.λ}
	\end{align}
\end{subequations} 
System \eqref{s.reduced} is called the \emph{core subsystem}. Consider the following ``error coordinates''
\begin{align}\label{d.error-coordinates}
	\tx\m &\coloneqq  x\m-\theta\sr, & \tx\p&=x\p, & \tz\p &\coloneqq  z\p-z\p\sr, & \tlam &\coloneqq  \lambda-\lambda\sr,
\end{align}
where $z\sr$ is defined in \eqref{d.z_sr}.
From  \eqref{e.dist_A_without_zm}, we obtain
\begin{equation}\label{e.sandwitch_dist_A}
	|(\tx\m,\tx\p,\tz\p,\tlam)| \le |(x,z,\lambda)|_{\cA\sr}\le \sqrt{n} |(\tx\m,\tx\p,\tz\p,\tlam)| \, .
\end{equation}
Hence, (semiglobal) exponential stability of the equilibrium point $(\theta\sr,0,z\p\sr,\lambda\sr)$ for \eqref{s.reduced} is equivalent to (semiglobal) exponential stability of $\cA\sr$ for Algorithm \eqref{s.algorithm}. 
Therefore, throughout the remainder of the proof, we restrict our attention to~\eqref{s.reduced}.

\subsection{Separation of the Optimization and the Network Dynamics}
This section highlights the presence of two interconnected dynamics in the core subsystem~\eqref{s.reduced} evolving at two different time scales. This sheds light on the internal structure of the core subsystem~\eqref{s.reduced} and lays the ground for the subsequent technical analysis.

Equations~\eqref{s.reduced.xm} and \eqref{s.reduced.λ} of the core subsystem~\eqref{s.reduced} can be written as
\begin{subequations}\label{s.opt-dyn}
	\begin{align}
		x\m^{t+1} &= x\m^t -  \gamma\1\T \nablaL(\1x\m^t,\lambda^t) + \gamma w^t \, ,
		\label{s.opt-dyn.x}
		\\
		\lambda^{t+1} &= \max\{0,\, \lambda^t + \gamma g(\1x\m^t) + \gamma d^t\} \, ,
		\label{s.opt-dyn.λ}
	\end{align}  
\end{subequations}
in which
	\begin{align}\label{d.w_d}
		w^t &\coloneqq   \1\T \big( \nablaL(\1x\m^t,\lambda^t)-\nablaL(x^t,\lambda^t)\big) \, ,&
		d^t &\coloneqq  g(x^t)-g(\1x\m^t) \, .
	\end{align} 
With $w^t=0$ and $d^t=0$, \eqref{s.opt-dyn} coincides with the \emph{centralized} Arrow-Hurwicz-Uzawa primal-dual algorithm analyzed in \cite{bin2024semiglobal}, which is proved to be semiglobally exponentially convergent to $(\theta\sr,\lambda\sr)$. System \eqref{s.opt-dyn} is therefore referred to as the \emph{optimization dynamics}.
In general, the presence of nonzero $w$ and $d$ leads to a perturbed version of the Arrow-Hurwicz-Uzawa algorithm that may be unstable. In view of \eqref{d.w_d}, within the core subsystem, the signals $w$ and $d$ are zero if  $x\p^t=0$ (namely, if consensus of the estimates $x_i^t$ has been reached). In other terms, Subsystems~\eqref{s.reduced.xm}, \eqref{s.reduced.λ} of the core subsystem~\eqref{s.reduced} can be seen as the (centralized) optimization dynamics~\eqref{s.opt-dyn} perturbed by the inputs~\eqref{d.w_d} that originate from the dynamics of the network seeking consensus.

Instead, the remaining subsystem \eqref{s.reduced.xp}-\eqref{s.reduced.zp} of the {core} subsystem \eqref{s.reduced} can be rewritten as  
\begin{equation}\label{s.cons-dyn}
	(x\p^{t+1},z\p^{t+1}) = A 	(x\p^{t},z\p^{t}) -\gamma {n}B \nablaL(\1\theta\sr,\lambda\sr) + \gamma {n} B v^t,
\end{equation}
in which
\begin{equation} \label{A_B}
	\begin{aligned}
		A & {\coloneqq} 
        \begin{bmatrix}
			I-S\T K S & -S\T K S\\
			S\T K S & I
		\end{bmatrix} \, ,  & 
        B & {\coloneqq} 
        \begin{bmatrix}
			S\T \\ 0
		\end{bmatrix} \, ,
	\end{aligned}
\end{equation}
and the input $v$ is defined as
\begin{equation}\label{d.v}
	v^t \coloneqq \nablaL(\1\theta\sr,\lambda\sr)-\nablaL(x^t,\lambda^t) \, .
\end{equation}
With $v^t=0$, Subsystem~\eqref{s.cons-dyn} describes the evolution of the dispersion component of the state variables $x$ and $z$; in particular, if $x\p^t \to 0$, the estimates $x_i^t$ reach consensus. Hence, \eqref{s.cons-dyn} is referred to as the \emph{network dynamics}.
As proved in \cite[Sec. IV-C]{binStabilityLinearConvergence2022}), $A$ is a {Schur} matrix; hence, if $v=0$, the dynamics~\eqref{s.cons-dyn} converges to an equilibrium. Indeed, a closer inspection reveals that $z\p$ acts as an integral action on $x\p$~\cite{binStabilityLinearConvergence2022,notarnicola_gradient_2023}, as 
\begin{equation*}
	z\p^{t+1} = z\p^t + S\T KS x\p^t \, ,
\end{equation*}
and $S\T K S$ is nonsingular in view of Assumption~\ref{ass.network}. It is then an immediate consequence that $x\p^t \to 0$. Namely, if $v=0$, $x$ reaches a consensual point.
Moreover, notice that $v^t =\nablaL(\1x\m^t,\lambda^t)-\nablaL( x^t,\lambda^t) +  \nablaL(\1\theta\sr,\lambda\sr)-\nablaL(\1x\m^t,\lambda^t)$, where the first difference vanishes with $x\p$ and, since $v^t$ is multiplied by $\gamma$ in \eqref{s.cons-dyn}, which is small, it can be dominated by the linear stable part. The second difference, instead, vanishes if the optimization dynamics reaches the optimum.
Thus, Subsystem~\eqref{s.reduced.xp}-\eqref{s.reduced.zp} can be seen as the network dynamics~\eqref{s.cons-dyn} perturbed by the input~\eqref{d.v} originating from the interconnection with the optimization dynamics~\eqref{s.opt-dyn}-\eqref{d.w_d}.

\begin{figure} 
	\centering
	\begin{tikzpicture}[font=\footnotesize,-latex']
		\tikzstyle{sys} = [draw,fill=black!5!white,line width=0.7pt,inner sep=1pt,minimum height=1.cm,text width=5cm,align=center] 
		\tikzstyle{line} = [draw, line width=0.7pt]
		\node[sys] (s1) at (0,0) {centralized optimization dynamics\\
			(\emph{$\gamma$-slow time scale})};
		
		\node[sys,below=.25cm of s1.south] (s2) {network dynamics\\(\emph{fast time scale})};

		\path[line] (s1.east) -| ($(s2.east)+(.5cm,0cm)$) -- (s2.east);		
		\path[line] (s2.west) -| ($(s1.west)-(.5cm,0cm)$) -- (s1.west); 
	\end{tikzpicture}
	\caption{Block scheme highlighting the time-scale separation of the optimization  and the network dynamics.}
	\label{fig.time-scales}
\end{figure}
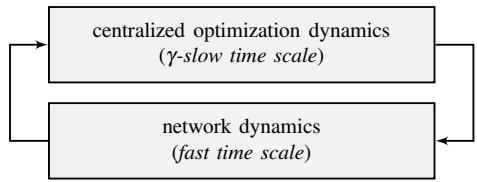

In summary, the core subsystem \eqref{s.reduced} can be seen as the interconnection of the perturbed optimization dynamics \eqref{s.opt-dyn}-\eqref{d.w_d} and the perturbed network dynamics  \eqref{s.cons-dyn}-\eqref{d.v} (see Figure~\ref{fig.time-scales}).
As each individual dynamics depends on the other though the perturbation terms, stability of their interconnection is not trivial, and the reminder of the proof is indeed aimed at establishing it.
The key observation on which the subsequent proof is based, is that the convergence rate of the optimization dynamics \eqref{s.opt-dyn} vanishes with $\gamma$.
Therefore, for small values of $\gamma$, the optimization dynamics evolves at a \emph{slow} time scale.  
Instead, $A$ does not depend on $\gamma$; hence, the convergence rate of the network dynamics \eqref{s.cons-dyn} is independent from $\gamma$. Thus, the network dynamics evolves at a \emph{fast} time scale. Moreover, $\gamma$ influences the \emph{gain} of both dynamics with respect to the other. These two key features enable an analysis approach based on the Lyapunov method and the small-gain principle.

We take on the analysis of the core subsystem by first considering the two perturbed dynamics separately; specifically, in Sections \ref{sec.proof.optimization_dynamic} and \ref{sec.proof.network_dynamics} we shall prove \emph{input-to-state stability (ISS)}~\cite{sontag1989smooth,jiang_input--state_2001} of the two dynamics with respect to the perturbation terms. 
Then, we put the two dynamics back together in Section~\ref{sec.proof.interconnection} thereby proving the main stability properties of the core subsystem.
Convergence and the exponential bound are then addressed in Section~\ref{sec.proof.convergence_rate}.
As a preliminary step, in the next Section~\ref{sec.proof.preliminary} we introduce some auxiliary quantities. 
In this respect, in the following we shall use the error coordinates \eqref{d.error-coordinates} along with the aggregate notation
\begin{align}\label{d.xi_txi}
	\xi&\coloneqq(x\p,z\p), &  \xi\sr &= (0, z\p\sr), & 	\txi&\coloneqq(\tx\p,\tz\p) = \xi-\xi\sr .
\end{align}

\subsection{Preliminary Definitions} \label{sec.proof.preliminary}
This section presents some preliminary definitions and results that will be used in the subsequent analysis. Grouping all such definitions here allows to more efficiently check the inter-dependencies between the various introduced quantities.
 
With $\Xi_0$ denoting the previously introduced compact initialization set, let $\kappa_0>0$ such that
\begin{equation}\label{d.kappa0-eq}
	\forall (x, z, \lambda) \in \Xi_0, \quad (x\m, \lambda, \xi) \in \overline{\mathbb{B}}_{\kappa_0} (\theta\sr, \lambda\sr, \xi\sr) \, .
\end{equation}
Moreover, let $P=P\T>0$ be the unique solution\footnote{Existence of $P$ follows since $A$ is Schur as shown in~\cite[Sec. IV-C]{binStabilityLinearConvergence2022}.} of the Lyapunov equation 
\begin{equation}\label{d.lyap-eq-P}
	A\T P A-P = -I \, ,
\end{equation}
and let  $c_l, c_u>0$ denote the smallest and largest eigenvalue of $P$, respectively.
With  
\begin{equation}\label{d.kappa1}
	{\kappa_1 \coloneqq \kappa_0 \sqrt{ \max\left\{3,\ \frac{c_u}{c_l} \right\}},}
\end{equation}
pick
\begin{equation}\label{d.kappa}
{\kappa \ge \max\left\{ 2\kappa_0+1 , \ \sqrt{ \frac{c_u}{c_l} } (  \kappa_1 + 1  ) \right\} },
\end{equation}
 and define
\begin{equation} \label{set_k}
	\cK  \coloneqq  \overline{\mathbb{B}}_{ \kappa}\left(\theta\sr, \lambda^{\star}, \xi\sr\right).
\end{equation}
Since $\mathcal{K}$ is compact, the smoothness properties guaranteed by  Assumption \ref{ass.optimization_problem_data} and  the optimality conditions \eqref{d.KKT} imply the existence of $k_0, k_1, k_2, k_3, k_4>0$  such that, for every   $(x\m,\lambda, \xi), (x\m',\lambda', \xi')\in\cK$, the following inequalities hold:   
\begin{equation} \label{d.k1_k2_k3_k4}
\begin{aligned}
{|\Phi(x,\lambda)-\Phi(x',\lambda')|} &\le k_0( |x\m-x\m'|+|x\p-x\p'| +|\lambda-\lambda'|) \, ,   \\ 
  |g(x)-g(x')| &\leq k_1|x-x'|  \, ,  \\
  {|g(\1x\m)-g(\1x\m')|} &\leq k_1|x\m-x\m'|  \, , \\
   |\1\T\nablaL(\1 x\m, \lambda)| &\leq k_2 \left( |x\m-\theta\sr|+|\lambda-\lambda^{\star}| \right) \, ,  \\
  { \left|\sum_{i\in\cN}\nabla f_i(x_i)-\sum_{i\in\cN}\nabla f_i(x_i')\right| }&\leq k_3|x - x'| \, ,  \\
  { \left|\sum_{i\in\cN}\nabla f_i(x\m)-\sum_{i\in\cN}\nabla f_i(x\m')\right| }&\leq k_3|x\m - x\m'| \, ,  \\
   |\nabla g(x)-\nabla g(x')| &\leq k_4|x - x'|    \, , 
   \\
   { |\nabla g(\1 x\m)-\nabla g(\1 x\m')|} &\leq k_4|x\m - x\m'| \, , 
\end{aligned}
\end{equation}
in which we used  
\begin{equation}\label{e.x_as_xm_xi}
	x=  \1 x\m + Sx\p,\qquad x\p= \begin{bmatrix}I&0\end{bmatrix}\xi,  
\end{equation}
and an analogous formula for $x'$ to express $x$ and $x'$ in terms of $(x\m,\xi)$ and $(x\m',\xi')$, respectively.
Moreover, we define the following constants
\begin{equation} \label{d.k5_k6_k7}
	\begin{aligned}
		k_5&\coloneqq    \sup_{(x\m,\lambda,\xi)\in\cK} \max \left\{|  \1\T\nablaL(x, \lambda)|, | \1\T \nablaL(\1 x\m, \lambda)| \right\} \, ,  
		\\
		 k_6&\coloneqq \sup_{(x\m,\lambda,\xi)\in\cK} \max \{ |\nabla g(x)|,   |\nabla g(\1 x\m)| \} \, ,  
		\\ 
		k_7&\coloneqq \sup_{(x\m,\lambda,\xi)\in\cK} \max \{ |g(x)|,   |g(\1 x\m)|\} \, , 
	\end{aligned}
\end{equation}
in which we have used again the previous expression of $x$ in terms of $x\m$ and $\xi$.

For every $i\in\cN$, let $m_i^a\le m_i$ denote the number of active constraints for agent $i$ at the optimum $x=\1\theta\sr$. Without loss of generality, we assume that these active constraints are associated with the  
indices $(i,j) \in \mathrm{\Lambda}_i \coloneqq  \{ (i,1), \ldots, (i,m_i^a) \}$. 
Thus, for every $i\in\cN$, we have $g_{i,j}\left(\theta\sr\right)=0$ for all $(i,j) \in \mathrm{\Lambda}_i$, and $g_{i,j}\left(\theta\sr\right)<0$ for all $(i,j) \in \mathrm{I}_i\coloneqq \left\{(i,m_i^a+1), \ldots, (i,m_i)\right\}$. 
Let $\mathrm{\Lambda}\coloneqq \cup_{i\in\cN} \mathrm{\Lambda}_i$ and $\mathrm{I}\coloneqq\cup_{i\in\cN} \mathrm{I}_i$, and let  $\lambda_{\mathrm{\Lambda}}\coloneqq \left(\lambda_{i,j}\right)_{(i,j)\in\mathrm{\Lambda}}$ collect all variables associated with active constraints and $\lambda_{\mathrm{I}}\coloneqq \left(\lambda_{i,j}\right)_{(i,j)\in\mathrm{I}}$ those associated with inactive constraints.
Moreover let $g_{\mathrm{\Lambda}}:\R^n\to\R^{m^a}$ and $g_{\mathrm{I}}:\R^n\to\R^{m-m^a}$, {with} $m^a\coloneqq m^a_1+\ldots+m^a_n$, be defined similarly as
\begin{align*}
	g_{\mathrm{\Lambda}}(x)&\coloneqq \left(g_{i,j}(x_i)\right)_{(i,j)\in\mathrm{\Lambda}},
	&
	g_{\mathrm{I}}(x)&\coloneqq \left(g_{i,j}(x_i)\right)_{(i,j)\in\mathrm{I}}.
\end{align*}
Then, \eqref{d.KKT.active_constraints} can be rephrased as
\begin{equation}
	g_{\mathrm{\Lambda}}\left(\1\theta\sr\right)=0, \quad \lambda_{\mathrm{\Lambda}}^{\star} \geq 0, \quad g_{\mathrm{I}}\left(\1\theta\sr\right)<0, \quad \lambda_{\mathrm{I}}^{\star}=0 .
\end{equation} 
Since $\nabla g$ is continuous then, according to Assumption \ref{ass.regular_θ*}, there exist $q_0 > 0$ and $\bar{\varepsilon}_1>0$ such that    
\begin{equation} \label{eq:q_0}
	\forall x \in \mathbb{R}^n,\quad {\left|x-\1\theta\sr\right|} \leq \bar{\varepsilon}_1 \implies  \begin{cases} 
		\nabla g_{\mathrm{\Lambda}}(x)^{\top} \nabla g_{\mathrm{\Lambda}}(x) \geq q_0 I, \\
		  g_{\mathrm{I}}(x)<0.
	\end{cases} 
\end{equation} 
Since {$\sum_{i\in\cN}f_i$} is strongly convex by Assumption \ref{ass.optimization_problem_data},
there exists $\mu_f>0$ such that 
\begin{equation}\label{d.mu_f}
	 \forall \theta \in \R, \qquad (\theta - \theta\sr) \left( \sum_{i\in\cN} \nabla f_i(\theta) - \sum_{i\in\cN} \nabla f_i(\theta\sr)  \right) \geq  \mu_f | \theta - \theta\sr |^2.
\end{equation} 
The constants defined below are adapted from \cite{bin2024semiglobal}, and are needed here to leverage the results of \cite{bin2024semiglobal} in the analysis of the optimization dynamics.
First, we fix
\begin{equation}\label{d.M_beta}
	\begin{aligned}
M &\ge \max \left\{ \frac{8}{k_2^2}\left( \frac{32}{3} c_u m^2(2k_1+1)^2+ \frac{9}{2c_u} k_0^2 n^2 \big(2|A\T PB|^2+|B\T PB|\big)\right), 1   \right\}  ,
\\
\beta&\coloneqq   \frac{{2(2+M)} k_2^2}{q_0},  
\end{aligned}
\end{equation}
where $q_0$ has been introduced in \eqref{eq:q_0}; then, we define the constants
\begin{equation} \label{d.alpha_delta}
\begin{aligned}
	\alpha_1  &\coloneqq \frac{q_0}{k_4 k_5+k_6\left(k_3+k_4\left|\lambda\sr\right|\right)},
	\quad  
	\alpha_2  \coloneqq \frac{k_4 k_5+k_6\left(k_3+k_4\left|\lambda^{\star}\right|\right)}{2 \alpha_1}+k_1 k_6, 
	\\
	\alpha_3  &\coloneqq \frac{3 k_1 k_2 k_6+k_2 k_4 k_5}{2}, 
	\quad
	\alpha_4 \coloneqq \frac{k_1 k_2 k_6+3 k_2 k_4 k_5}{2}, 
	\\
	\alpha_7   &\coloneqq k_2 k_6 + \frac{\beta}{2}\left( \kappa k_4 k_6+k_4 k_5+k_2 k_6\right), 
	\quad
	\alpha_8   \coloneqq \frac{\beta}{2}\left(k_1^2 k_6^2+ \kappa k_4 k_6 k_2^2\right), 
	\\ 
	\alpha_9 &\coloneqq \beta\left(\frac{k_2 k_6}{2}+k_6^2\right)+k_2 k_6, 
	\quad
	\alpha_{10}\coloneqq \beta \frac{ \kappa k_4 k_6 k_2^2}{2},
	\\
	\delta_1&\coloneqq \frac{M k_2^2}{8 \alpha_9},  \quad {\delta_2\coloneqq  \begin{cases}0 & \text { if } m=m^a \\ \frac{{\mu_f}}{2\left(m-m^a\right)(1+\beta) k_1} & \text { if } m> m^a \end{cases} }
	\\
	\alpha_{11}  &\coloneqq \frac{\beta}{2}\left(k_4 k_5+k_2 k_6 \frac{1+\delta_1+\delta_1^2}{\delta_1}\right)+k_2 k_6 \frac{1+\delta_1}{\delta_1} 
	+k_6^2+\beta k_6^2 \frac{1+2 \delta_1}{4 \delta_1}+\beta \kappa k_4 k_6.
\end{aligned}
\end{equation}   
Next, with
\begin{align*} 
	h&\coloneqq \min _{(i,j) \in \mathrm{I}}\left|g_{i,j}\left( {\theta\sr}\right)\right|, 
	&
	\bar{\varepsilon}&\coloneqq \min \left\{\bar{\varepsilon}_1, \frac{h}{4 k_1(1+2 \beta)}\right\} ,
\end{align*} 
let us arbitrarily pick
\begin{equation}\label{d.epsilon}
	\varepsilon \in(0, \bar{\varepsilon}),
\end{equation}
and define the constants 
\begin{equation}\label{d.bargamma}
\begin{aligned}
	\bar{\gamma}_1&\coloneqq  \frac{1}{16 \kappa_0 k_5}  , \quad \bar{\gamma}_2\coloneqq \frac{1}{2 k_5}, 
	 \quad
	 \bar{\gamma}_3\coloneqq \frac{1}{16 \kappa_0 k_7}, 
	 \quad 
	 \bar{\gamma}_4\coloneqq {\frac{1}{2 k_7}}, 
	\\
	\bar{\gamma}_5&\coloneqq \frac{1}{\beta k_6}, 
	\quad 
	\bar{\gamma}_6\coloneqq \frac{\mu_f \varepsilon^2}{\beta\left(4 \kappa_0^2 k_4 k_5+2 \kappa_0 k_6 k_7+k_5 k_6 \kappa \right)+k_7^2+k_5^2},
	\\
	\bar{\gamma}_7 &\coloneqq \min \left\{\frac{\mu_f}{2\left(k_1^2+2 k_2^2+\beta \alpha_2\right)}, \sqrt{\frac{\mu_f}{2 \beta \alpha_3}}\right\}, 
	\quad 
	\bar{\gamma}_8\coloneqq \frac{M k_2^2}{2 \beta \alpha_4}, \\
	\bar{\gamma}_9&\coloneqq \min \left\{\frac{1}{2 \sqrt{\alpha_{11}}}, \frac{\delta_2}{4(1+\beta) k_1}\right\}, 
	\quad \bar{\gamma}_{10}\coloneqq \frac{h}{4 \alpha_{11} \kappa}, \\
	\bar{\gamma}_{11}&\coloneqq \min \left\{\frac{\mu_f}{8 \alpha_7}, \frac{1}{2} \sqrt[3]{\frac{\mu_f}{{\alpha_{8}}}},  {\sqrt{\frac{Mk_2^2}{8\alpha_{10}}}} \right\}, 
	\quad \bar{\gamma}_{12}\coloneqq \frac{2h}{\kappa k_2^2}.
\end{aligned}
\end{equation}
In addition to these constants adapted from~\cite{bin2024semiglobal}, we introduce the following ones (recall that $A$ and $B$ are defined in \eqref{A_B})\footnote{The reader can readily verify that the above-defined constants contain no circular definitions.}
\begin{equation}\label{d.bargamma_2}
	\begin{aligned} 
		 \bar{\gamma}_{13}&\coloneqq \min\left\{ \frac{1}{12\kappa_1 k_5}, \ \frac{1}{\sqrt{6}k_5},\  \frac{1}{12\kappa_1 k_7}, \ \frac{1}{\sqrt{6}k_7} \right\} \, , \\  \bar\gamma_{14}  &\coloneqq  \frac{ \min\left\{\sqrt{\frac{c_l}{3}}, \frac{1}{2}  \right\}}{nk_0 {\kappa} \sqrt{6|A\T PB|^2+3|B\T PB|} }  \, , \\ 
		 \bar\gamma_{15} &\coloneqq \min\left\{ \frac{1}{2\beta k_0(n\kappa_1 {k_4}+\sqrt{n}{k_6})} ,\, \frac{1}{\sqrt{2\beta nk_0k_4k_5}} \right\},\\
		 \bar\gamma_{16}&\coloneqq\min\left\{ \frac{1}{{2}k_2}, \frac{\sqrt{n}{k_0}}{\beta k_1k_6+mk_1(2k_1+1)},\ \sqrt{\frac{2\sqrt{n} k_0}{\beta k_1k_2k_6}},\, \frac{m(2k_1+1)}{4\sqrt{n}k_0k_2},\,\sqrt{\frac{m(2k_1+1)}{2\beta k_1k_2k_6}} \right\},
		 \\
		 \bar\gamma_{17}&\coloneqq\min\left\{ \frac{\mu_f}{768 nk_0^2(2c_u+3)},\, \sqrt{\frac{3}{64 c_u (nk_0^2+k_1^2+ k_1)}} \right\},
		 \\
		 \bar\gamma_{18} &\coloneqq \frac{\mu_f c_u}{12(2c_u+3)k_0^2n^2(2|A\T PB|^2 + |B\T PB|)},
		 \quad 
		 \bar\gamma_{19}\coloneqq  
		\frac{6\mu_f \varepsilon^2}{Mk_2^2 \kappa_1^2} ,
		\\
		\bar\gamma_{20}&\coloneqq \min\left\{\sqrt{\frac{2}{Mk_2^2}},\, \frac{2h}{Mk_2^2\kappa_1} \right\} \, , 
	\end{aligned}
\end{equation} 
and we set
\begin{equation} \label{gamma_all}
	\bar\gamma_0\coloneqq \min _{i=1, \ldots, 20} \bar\gamma_i.
\end{equation}
Finally, we fix once and for all 
\begin{equation}\label{d.gamma}
	\gamma \in (0, {\bar\gamma}).
\end{equation}
Now that all these auxiliary constants are defined, we start analyzing the core subsystem~\eqref{s.reduced}. In particular, Sections~\ref{sec.proof.optimization_dynamic} and \ref{sec.proof.network_dynamics} analyze the perturbed optimization and network dynamics, respectively.  Section~\ref{sec.proof.interconnection} combines the obtained results thereby proving some stability properties of the core subsystem. Finally, Section~\ref{sec.proof.convergence_rate} refines such a result by providing the sought exponential bound. 

In the remainder of the proof, we shall often refer to the above-defined constants.
Moreover, we will make use of the Young's inequality
\begin{equation}\label{d.Young}
	2ab \le \frac{a^2}{\epsilon} + \epsilon b^2 
\end{equation}
that holds for all $a,b\in\R$ and $\epsilon>0$.

\subsection{The Optimization Dynamics} \label{sec.proof.optimization_dynamic}
In this section, we analyze the optimization dynamics~\eqref{s.opt-dyn}-\eqref{d.w_d}. We show that such a dynamics is semiglobally ISS relative to the optimal equilibrium $(\theta\sr,\lambda\sr)$ and with respect to the perturbation terms $w$ and $d$, which are treated here as exogenous variables.
To simplify the notation, we use the error coordinates \eqref{d.error-coordinates}, \eqref{d.xi_txi}.

The starting point is the result of \cite{bin2024semiglobal} that establishes semiglobal exponential stability of the \emph{unperturbed} (i.e., with $w^t=0$ and $d^t=0$) optimization dynamics \eqref{s.opt-dyn}-\eqref{d.w_d}. The latter can be compactly written as
\begin{equation}\label{s.opt-dyn-unp}
	(x\m^{t+1},\lambda^{t+1}) = \phi_\gamma( x\m^t,\lambda^t),
\end{equation}
where 
\begin{align*}
	\phi_\gamma( x\m ,\lambda) &\coloneqq  \left( x\m  -   \gamma\1\T \nablaL(\1x\m ,\lambda ) ,\  \max\{0,\, \lambda  + \gamma g(\1x\m)  \} \right ). 
\end{align*}   
To prove semiglobal exponential stability of~\eqref{s.opt-dyn-unp},~\cite{bin2024semiglobal} proposed the following Lyapunov function
\begin{equation}\label{d.V-opt}
		V_{\rm opt}(x\m,\lambda)  \coloneqq   |x\m-\theta\sr|^2 + |\lambda-\lambda\sr|^2 
		 + \gamma\beta (x\m-\theta\sr) {\nabla g(\1 x\m)} (\lambda-\lambda\sr) \, ,
\end{equation}
for a suitable value of $\beta$ here taken as in~\eqref{d.M_beta}.
By defining the level set
\begin{equation}\label{d.Omega_opt}
	\Omega_{\rm opt} \coloneqq \left\{ (x\m,\lambda)\in\overline{\mathbb{B}}_{\kappa}(\theta\sr,\lambda\sr)\st V_{\rm opt}(x\m,\lambda) \le \frac{3}{2}\kappa_0^2 \right \} \, ,
\end{equation}
the following result can be deduced from~\cite{bin2024semiglobal}.

\begin{lemma}\label{lem.MP}
Suppose that Assumptions~\ref{ass.optimization_problem_data}-\ref{ass.regular_θ*} hold and let $\varepsilon$ be given by~\eqref{d.epsilon}. Then: 
\begin{enumerate}
	\item\label{item.lem.Mp.sandwitch} For every $(x\m,\lambda)\in\overline{\mathbb{B}}_{\kappa}(\theta\sr,\lambda\sr)$,
	\begin{equation}\label{e.sandwich-Vopt}
		\frac{1}{2}|(x\m-\theta\sr,\lambda-\lambda\sr)|^2 \le V_{\rm opt}(x\m,\lambda)\le \frac{3}{2} |(x\m-\theta\sr,\lambda-\lambda\sr)|^2.
	\end{equation}  
	\item\label{item.lem.Mp.ge_epsilon} For every $(x\m,\lambda)\in 	\Omega_{\rm opt}$ such that $|x\m-\theta\sr|>\varepsilon$, 
	\begin{equation} \label{e.Vopt_ge_epsilon}
		V_{\rm opt}\big(\phi_\gamma(x\m,\lambda)\big)\le  V_{\rm opt}(x\m,\lambda) {-2\gamma\mu_f|x\m-\theta\sr|^2+\gamma \mu_f \varepsilon^2.}
	\end{equation} 
	
	\item\label{item.lem.Mp.le_epsilon} For every $(x\m,\lambda)\in 	\Omega_{\rm opt}$ such that $|x\m-\theta\sr|\le\varepsilon$,
	\begin{equation}
			\begin{aligned}
			 V_{\rm opt}\big(\phi_\gamma(x\m,\lambda)\big)&
            \le 
            V_{\rm opt}(x\m,\lambda)  -\frac{1}{2} \sum_{i \in \mathrm{I}} \min \left\{\left|\lambda_i\right|^2, \gamma h\left|\lambda_i\right|\right\}
			 \\&\qquad 
		 -\frac{1}{4} \gamma \mu_f\left|x_{\mathrm{m}}-\theta^{\star}\right|^2-\frac{1}{4} \gamma^2  M k_2^2\left|\lambda_{\mathrm{\Lambda}}-\lambda_{\mathrm{\Lambda}}^{\star}\right|^2.
		\end{aligned}
	\end{equation}
\end{enumerate} 
\end{lemma}
\begin{proof} 
	System \eqref{s.opt-dyn-unp} coincides with the algorithm analyzed in \cite{bin2024semiglobal} with $f\coloneqq\sum_{i\in\cN}f_i$, with the same function $g$, and with \cite[Ass. 1]{bin2024semiglobal} and \cite[Ass. 2]{bin2024semiglobal} that are implied by Assumptions~\ref{ass.optimization_problem_data} and \ref{ass.regular_θ*}.
	Regarding Items \ref{item.lem.Mp.ge_epsilon} and \ref{item.lem.Mp.le_epsilon}, we discuss the two cases $M=1$ and $M>1$.
	For $M=1$, the condition $\gamma<\bar\gamma_0$ implies $\gamma<\bar\gamma$ (i.e., Equation (20)) in \cite{bin2024semiglobal} since, for $M=1$, all constants defined in previous Section~\ref{sec.proof.preliminary} are upper bounds of those defined in \cite[Sec.~4.1]{bin2024semiglobal} (the only different notations are $K_0,K,q,c_0,r,r_a$, and $\Omega_\rho$ in \cite{bin2024semiglobal} that must be set as $K_0\coloneqq\kappa_0$, $K\coloneqq\kappa$, $q\coloneqq q_0$, $c_0\coloneqq\mu_f$, $r\coloneqq m$, $r_a \coloneqq m^a$, and $\Omega_\rho\coloneqq\Omega_{\rm opt}$).
	Thus, for $M=1$, the claim of the lemma directly follows from~\cite{bin2024semiglobal}; specifically, see the inequality before \cite[Eq. (34)]{bin2024semiglobal} in \cite[Sec.~4.4 ]{bin2024semiglobal} and  \cite[Eq. (61)]{bin2024semiglobal} in \cite[Sec.~4.5 ]{bin2024semiglobal}.
	
	Instead, for any $M>1$, one can obtain the claim of the lemma by means of the same arguments of \cite{bin2024semiglobal}, by redefining $\beta$, $\delta_1$, $\bar\gamma_8$, and {$\bar\gamma_{11}$} as done in Section~\ref{sec.proof.preliminary} above.
	Indeed, if $\beta=\frac{2}{q_0}(2+M)k_2^2$ and $\gamma<\bar\gamma_8=\frac{Mk_2^2}{2\beta\alpha_4}$ in \cite{bin2024semiglobal}, one gets by \cite[Eq. (47)]{bin2024semiglobal} the term $-\gamma^2 \frac{M k_2^2}{2} |\tlam_A|^2$ in place of  $-\gamma^2 \frac{k_2^2}{2} |\tlam_A|^2$. The same substitution is inherited by \cite[Eq. (54)]{bin2024semiglobal}, \cite[Eq. (57)]{bin2024semiglobal}, and then \cite[Eq. (58)]{bin2024semiglobal}. From such an inequality, it is then clear that, if $\delta_1\alpha_9 = \frac{M}{8}k_2^2$ (which is implied by the new definition of $\delta_1$ in \eqref{d.alpha_delta}) and $\gamma^2\alpha_{10}\le \frac{M}{8}k_2^2$ (which is implied by $\gamma<\bar\gamma_{11}$ as defined in \eqref{d.bargamma}, and in particular by $\bar\gamma_{11}\le (\frac{M}{8\alpha_{10}}k_2^2)^{1/2}$), then the term multiplying $|\tlam_A|^2$ in \cite[Eq. (58)]{bin2024semiglobal} can be upper-bounded by $-\frac{\gamma^2Mk_2^2}{4}|\tlam_A|^2$.
    The proof follows since the passages after \cite[Eq. (58)]{bin2024semiglobal} (until \cite[Eq. (61)]{bin2024semiglobal}) do not introduce other terms dependent from $\tlam_A$. Finally, in both cases   Item \ref{item.lem.Mp.sandwitch} follows from  \cite[Lem.~6]{bin2024semiglobal} that is not affected by $M$.
    \end{proof}

Next, we consider the perturbed optimization dynamics \eqref{s.opt-dyn} that can be compactly rewritten as
\begin{equation}\label{s.opt-dyn-perturbed}
	(x\m^{t+1},\lambda^{t+1})  = \tilde\phi_\gamma(x\m^t,\lambda^t,w^t,d^t),   	
\end{equation} 
in which
\begin{equation*}
\tilde\phi_\gamma(x\m,\lambda,w,d)  	\coloneqq 
  \left(
		x\m  -  \gamma\1\T \nablaL(\1x\m ,\lambda ) +\gamma w,\ 
		\max\{0,\, \lambda  + \gamma g(\1x\m)+ \gamma d  \}
	\right)     ,
\end{equation*}
and where $w$ and $d$ are treated here as exogenous inputs.
%
%
The $1$-step update of the Lyapunov candidate \eqref{d.V-opt} along the solutions of the perturbed optimization dynamics~\eqref{s.opt-dyn-perturbed} {can be rewritten as}
\begin{equation}\label{e.Delta_Vopt_tilde}
	V_{\rm opt}\big(\tilde\phi_\gamma(x\m,\lambda,w,d)\big)  = V_{\rm opt}\big(\phi_\gamma(x\m,\lambda)\big)  +  W_\gamma(x\m,\lambda,w,d) ,
\end{equation}
where 
\begin{equation}\label{d.W}
	W_\gamma(x\m,\lambda,w,d)\coloneqq  V_{\rm opt}\big(\tilde\phi_\gamma(x\m,\lambda,w,d)\big)-V_{\rm opt}\big(\phi_\gamma(x\m,\lambda)\big)
\end{equation}
equals the mismatch between the unperturbed and the perturbed updates. 
Since $W_\gamma$ is continuous and vanishes with $(w,d)$, it is then an immediate consequence of Lemma~\ref{lem.MP} that the perturbed optimization dynamics \eqref{s.opt-dyn-perturbed} is locally ISS with respect to $w$ and $d$.
In addition, it can be shown that it is semiglobally ISS, i.e., that if $w$ and $d$ are bounded, then for each compact initialization set, and for all small-enough values of $\gamma$, \eqref{s.opt-dyn-perturbed} is ISS with a restriction on the initial conditions that need to lie within the specified initialization set. This latter property will be implicitly used in Section~\ref{sec.proof.interconnection} to deal with the interconnection between the optimization and the network dynamics.   

 \subsection{The Network Dynamics}\label{sec.proof.network_dynamics}
  This section considers the network dynamics~\eqref{s.cons-dyn}; it is shown that such a system is globally ISS relative to the target equilibrium $\xi\sr$ and with respect to the input $v$ (we recall that $\xi\sr$ is defined in \eqref{d.xi_txi} with $z\p\sr$ given by \eqref{d.z_sr}).
 Consider the Lyapunov candidate
 \begin{equation}
 	V_{\rm net}(\xi) \coloneqq (\xi-\xi\sr)\T P(\xi-\xi\sr) \, ,
 \end{equation}
 where, we recall, $P$ is defined in Section~\ref{sec.proof.preliminary} as the unique symmetric and positive-definite solution of~\eqref{d.lyap-eq-P}. 
 By treating $v$ as a generic exogenous input, and using the fact that, by construction
 \begin{equation*}
 	\xi\sr + \gamma n B\Phi(\1\theta\sr,\lambda\sr) = A\xi\sr,
 \end{equation*} 
one obtains that the increment of $V_{\rm net}$ along the solutions of \eqref{s.cons-dyn} satisfies 
  	\begin{align*}
 		V_{\rm net}&\big(A\xi - \gamma n B \Phi(\1 \theta\sr,\lambda\sr) +\gamma nB v\big) -V_{\rm net}(\xi)  
\\& = (A(\xi-\xi\sr) +\gamma nB v)\T P(A(\xi-\xi\sr) +\gamma nB v)  - (\xi-\xi\sr)\T P(\xi-\xi\sr)
 		\\
 		&= - |\txi|^2 + 2\gamma n \txi\T A\T P B v + \gamma^2n^2 v\T B\T P B v \, .
 	\end{align*}
 By using
 \begin{equation}\label{e.bounds_P}
 \forall \xi\in\R^{2(n-1)},\quad 	c_l |\xi-\xi\sr|^2 \le V_{\rm net}(\xi) \le c_u |\xi-\xi\sr|^2,
 \end{equation}
 where, we recall, $c_l$ and $c_u$ have been introduced in Section~\ref{sec.proof.preliminary} as the smallest and largest eigenvalues of $P$, respectively,  and using the Young's inequality \eqref{d.Young} with $a=|\txi|$, $b=\gamma n|A\T P B||v|$, and $\epsilon=2$, we obtain
 \begin{equation}\label{e.V_net_1}
 	\begin{aligned}
 		V_{\rm net}&\big(A\xi - \gamma n B \Phi(\1 \theta\sr,\lambda\sr) +\gamma nB v\big) -V_{\rm net}(\xi)  
 		\\
 		&\qquad\le  - |\txi|^2 +  2n\gamma  |\txi|  |A\T P B| |v| + \gamma^2n^2 |B\T P B| |v|^2\\
 		&\qquad\le -\frac{1}{2} |\txi|^2 + \gamma^2n^2\left(2|A\T P B|^2 +|B\T P B|\right) |v|^2,
 	\end{aligned}  
 \end{equation}
 from which we can deduce that the perturbed network dynamics \eqref{s.cons-dyn} is ISS relative to the target equilibrium $\xi\sr$ and with respect to the input $v$.
 
 \subsection{Interconnection Between the Optimization and Network Dynamics} \label{sec.proof.interconnection}
This section considers the interconnection between the optimization and the network dynamics obtained by letting $w$, $d$, and $v$ be given by Equations~\eqref{d.w_d} and \eqref{d.v}, respectively. For simplicity, in the following we use the notation
\begin{align}\label{e.notation_plus}
	{x\m^+} &\coloneqq x\m - \gamma \1\T \Phi(x,\lambda), & \xi^+&\coloneqq A\xi -\gamma n B \Phi(x,\lambda) , &
	\lambda^+ &\coloneqq \max\{ 0,\ \lambda+\gamma g(x)\}
\end{align}
to denote the update map of the core subsystem~\eqref{s.reduced}.
We analyze this interconnection through a candidate Lyapunov function consisting of a weighted sum of the optimization candidate $V_{\rm opt}$ and the network candidate $V_{\rm net}$, i.e.,
\begin{equation} \label{d.V}
	V(x\m, \lambda, \xi) \coloneqq V_{\rm opt}(x\m, \lambda)+ {\frac{3}{2c_u}}V_{\rm net}(\xi).
\end{equation}  
Define the level set
\begin{equation}
	\Omega \coloneqq \left\{(x\m, \lambda,\xi) \in\cK: V(x\m, \lambda, \xi)  \leq   { \frac{3}{2} } \kappa_0^2 \right\},
\end{equation}
then, the following result holds. 
\begin{lemma}\label{lem.Omega}
	Suppose that Assumptions~\ref{ass.optimization_problem_data},  \ref{ass.regular_θ*}, and \ref{ass.network} hold. Then:
	\begin{enumerate}
		\item\label{item.lem.Omega.sandwich_Omega} $\overline{\mathbb{B}}_{\kappa_0}\left(\theta^{\star}, \lambda^{\star}, \xi\sr\right) \subset \Omega  \subset \overline{\mathbb{B}}_{ \kappa_1}\left(\theta^{\star}, \lambda^{\star}, \xi\sr\right)$.
		\item\label{item.lem.Omega.plus}  $\forall (x\m, \lambda, \xi ) \in \overline{\mathbb{B}}_{\kappa_1 }(\theta^{\star}, \lambda^{\star}, \xi\sr )$, $(x\m^{+}, \lambda^{+}, \xi^+  ) \in \cK$.
		\item \label{item.lem.Omega.opt}  $\forall(x\m, \lambda, \xi ) \in\Omega$, $(x\m,\lambda)\in\Omega_{\rm opt}$.
	\end{enumerate}    
\end{lemma}
\begin{proof}
    See Appendix~\ref{proof.lem.Omega}. 
\end{proof}
Lemma~\ref{lem.Omega} implies that, for every initial condition  of Algorithm \eqref{s.algorithm} belonging to $\Xi_0$, the corresponding initial condition of the core subsystem lies inside $\Omega$; in turn, for every $(x\m,\lambda,\xi)\in\Omega$, the bounds in~\eqref{d.k1_k2_k3_k4} and~\eqref{d.k5_k6_k7} established in Section~\ref{sec.proof.preliminary} hold. Moreover, the same bounds also hold for $(x\m^+, \lambda^+, \xi^+)$ as it lies in~$\cK$.

Next, we analyze the increment of $V$ along the solutions of the core subsystem. 
Starting from the preliminary computations made in Sections~\ref{sec.proof.optimization_dynamic} and~\ref{sec.proof.network_dynamics}, we first analyze separately the increments of $V_{\rm opt}$ and $V_{\rm net}$ when $(w,d)$ and $v$ are given by~\eqref{d.w_d} and~\eqref{d.v}, respectively, and then we combine the results to obtain a bound for the increment of $V$.
We start from the optimization dynamics and, in particular, from~\eqref{d.W} when evaluated at $(w,d)$ as given in~\eqref{d.w_d}. Formally, let 
\begin{equation}\label{d.bar_W}
	\overline{W}_\gamma (x\m,\lambda,\xi)\coloneqq W_\gamma\left(x\m,\lambda,\1\T \big( \nablaL(\1x\m,\lambda)-\nablaL(x,\lambda)\big), g(x)-g(\1x\m) \right)  \, ,
\end{equation} 
%
where we implicitly use~\eqref{e.x_as_xm_xi} to express $x$ in terms of $x\m$ and $\xi$. 
The following lemma provides a bound $\overline{W}_\gamma (x\m,\lambda,\xi)$ for every
$(x\m,\lambda,\xi)\in\overline{\mathbb{B}}_{\kappa_1}(\theta\sr,\lambda\sr,\xi\sr)$.
\begin{lemma} \label{lem.W}
	Suppose that  Assumptions~\ref{ass.optimization_problem_data}, \ref{ass.regular_θ*}, and \ref{ass.network} hold. Then, for every $(x\m,\lambda,\xi)\in\overline{\mathbb{B}}_{\kappa_1}(\theta\sr,\lambda\sr,\xi\sr)$, the following estimate holds:
		\begin{equation}  \label{eq.w}
			\overline{W}_\gamma (x\m,\lambda,\xi)  
			 \leq  \gamma \frac{\mu_f c_u}{4(2c_u+3)}|\tx\m|^2  +  \frac{32}{3}\gamma^2 c_u m^2(2k_1+1)^2 |\tlam|^2 + \frac{3}{16 c_u}|\txi|^2 \, . 
	\end{equation}
\end{lemma}
\begin{proof}
    See Appendix~\ref{proof.lem.W}.
\end{proof}
 
 Plugging~\eqref{d.w_d} into  \eqref{s.opt-dyn-perturbed} and \eqref{e.Delta_Vopt_tilde}, and using~\eqref{eq.w}, lead to the following bound on the update of $V_{\rm opt}$  
 \begin{equation}\label{e.Delta_Vopt_tilde_2}
 	\begin{aligned}
 		V_{\rm opt}\big(x\m^+,\lambda^+,\xi^+\big)  &\le  V_{\rm opt}\big(\phi_\gamma(x\m,\lambda)\big)  +  \gamma \frac{\mu_f c_u}{4(2c_u+3)}|\tx\m|^2 \\
 		&\qquad +  \frac{32}{3}\gamma^2 c_u m^2(2k_1+1)^2 |\tlam|^2 + \frac{3}{16 c_u}|\txi|^2 \, ,
 	\end{aligned}
 \end{equation}
 for every
 $(x\m,\lambda,\xi)\in\overline{\mathbb{B}}_{\kappa_1}(\theta\sr,\lambda\sr,\xi\sr)$. 
 Note that the first term in the right-hand side of~\eqref{e.Delta_Vopt_tilde_2} is the update of $V_{\rm opt}$ in the centralized case, without network. In turn, Bound~\eqref{e.Delta_Vopt_tilde_2} estimates the destabilizing effect on the optimization dynamics due to the network dynamics.

As far as the network dynamics is concerned, we proceed as in the proof of Lemma~\ref{lem.Omega}. By substituting $v$ with \eqref{d.v} in \eqref{e.V_net_1}, and performing the same computations as in Inequality~\eqref{e.proof.Omega.Vnet} in the proof of Lemma~\ref{lem.Omega}, we obtain 
%
%
\begin{equation*}
	V_{\rm net}(\xi^+)  \le V_{\rm net}(\xi) -\frac{1}{4}|\txi|^2+   3\gamma^2 k_0^2 n^2 \big(2|A\T PB|^2+|B\T PB|\big) \big( |\tx\m|^2  + |\tlam|^2\big) \, ,
\end{equation*}
for every
$(x\m,\lambda,\xi)\in\overline{\mathbb{B}}_{\kappa_1}(\theta\sr,\lambda\sr,\xi\sr)$. Using $\gamma<\bar\gamma_{18}$ yields
\begin{equation}\label{e.Vnet_last}
	\begin{aligned}
			V_{\rm net}(\xi^+)  &\le V_{\rm net}(\xi) -\frac{1}{4}|\txi|^2+   \gamma \frac{\mu_f c_u}{4(2c_u+3)} |\tx\m|^2 \\
			&\qquad+ 3\gamma^2 k_0^2 n^2 \big(2|A\T PB|^2+|B\T PB|\big)   |\tlam|^2 \, ,
	\end{aligned} 
\end{equation}
for every
$(x\m,\lambda,\xi)\in\overline{\mathbb{B}}_{\kappa_1}(\theta\sr,\lambda\sr,\xi\sr)$.
Inequality~\eqref{e.Vnet_last} estimates the destabilizing effect of the optimization dynamics on the network dynamics.

Inequalities \eqref{e.Delta_Vopt_tilde_2} and  \eqref{e.Vnet_last} are now combined, and it is shown that, when combined, the destabilization effects each dynamics yields on the other balance out, thereby leading to a stable interconnection.
By combining \eqref{e.Delta_Vopt_tilde_2} and  \eqref{e.Vnet_last}, and using the definition of $M$ in \eqref{d.M_beta}, from \eqref{d.V} one obtains
\begin{equation}\label{e.V_1}
	\begin{aligned}
		V(x\m^+,\lambda^+,\xi^+) &\le  V_{\rm opt}\big(\phi_\gamma(x\m,\lambda)\big) + \frac{3}{2c_u}V_{\rm net}(\xi)\\
		&\qquad +  \gamma   \left(1+ \frac{3}{2c_u}\right)\frac{\mu_f c_u}{4(2c_u+3)}  |\tx\m|^2   - \frac{3}{16 c_u}|\txi|^2  
		\\
		&\qquad + \gamma^2 \left( \frac{32}{3} c_u m^2(2k_1+1)^2+ \frac{9}{2c_u} k_0^2 n^2 \big(2|A\T PB|^2+|B\T PB|\big)\right) |\tlam|^2 
		\\
		&\le V_{\rm opt}\big(\phi_\gamma(x\m,\lambda)\big) + \frac{3}{2c_u}V_{\rm net}(\xi)
		\\
		&\qquad +  \gamma  \frac{\mu_f}{8}  |\tx\m|^2   - \frac{3}{16 c_u}|\txi|^2 + \frac{1}{8}\gamma^2  M k_2^2 |\tlam|^2,
	\end{aligned} 
\end{equation}
for every
$(x\m,\lambda,\xi)\in\overline{\mathbb{B}}_{\kappa_1}(\theta\sr,\lambda\sr,\xi\sr)$.

The term $V_{\rm opt}\big(\phi_\gamma(x\m,\lambda)\big)$ in the right-hand side of the previous bound is the update of the optimization Lyapunov candidate evaluated along the solutions of the centralized optimization dynamics. We now bound such a term by leveraging the results of Section~\ref{sec.proof.optimization_dynamic} (in particular, Lemma~\ref{lem.MP}). Pick $(x\m,\lambda,\xi)\in\Omega$; in view of Lemma~\ref{lem.Omega},  $\Omega\subset
\overline{\mathbb{B}}_{\kappa_1}(\theta\sr,\lambda\sr,\xi\sr)$, so as all the previous bounds apply to $(x\m,\lambda,\xi)$.
Moreover, Lemma~\ref{lem.Omega} also implies $(x\m,\lambda)\in\Omega_{\rm opt}$.
We now analyze separately the cases $|x\m-\theta\sr|> \varepsilon$ and $|x\m-\theta\sr|\le \varepsilon$.

Firstly, suppose that $|x\m-\theta\sr|> \varepsilon$.
Then, $-|\tx\m|^2\le -\varepsilon^2$; hence, by using Item \ref{item.lem.Mp.ge_epsilon}  of Lemma~\ref{lem.MP} and $\gamma<\bar\gamma_{19}$, and since $|\tlam|\le \kappa_1$, we get:
\begin{equation}\label{e.V_g_epsilon}
	 \begin{aligned}
		V(x\m^+,\lambda^+,\xi^+)  
		&\le V_{\rm opt}(x\m,\lambda)  -2\gamma\mu_f|\tx\m|^2+\gamma \mu_f \varepsilon^2  + \frac{3}{2c_u}V_{\rm net}(\xi)
		\\
		&\qquad +  \gamma  \frac{\mu_f}{8}  |\tx\m|^2   - \frac{3}{16 c_u}|\txi|^2 + \frac{1}{8} \gamma^2  M k_2^2 |\tlam|^2,
		\\
		&\le V(x\m,\lambda,\xi) - \frac{15}{8} \gamma \mu_f |\tx\m|^2 +\gamma\mu_f\varepsilon^2 +\frac{1}{8}\gamma^2 Mk_2^2 \kappa_1^2
		\\
		&\le V(x\m,\lambda,\xi) - \frac{1}{8}\gamma\mu_f \varepsilon^2 \, .
	\end{aligned}
\end{equation}

 Next, suppose that $|x\m-\theta\sr|\le\varepsilon$. Then, from~\eqref{e.V_1} and Item~\ref{item.lem.Mp.le_epsilon} of Lemma~\ref{lem.MP}, it follows that
 \begin{equation}\label{e.V_le_eps}
 	\begin{aligned}
    V(x\m^+,\lambda^+,\xi^+)  	&\le  V(x\m,\lambda,\xi)  
    -  \gamma  \frac{\mu_f}{8}  |\tx\m|^2 
    -\frac{1}{8} \gamma^2 Mk_2^2|\tlam_{\mathrm{\Lambda}}|^2 
    - \frac{3}{16 c_u}|\txi|^2
 		\\
 		&\qquad 
 		+ \frac{1}{8}\gamma^2  M k_2^2 |\tlam_{\mathrm{I}}|^2 -\frac{1}{2} \sum_{i \in {\mathrm{I}}} \min \left\{|\tlam_i|^2, \gamma h|\tlam_i|\right\} .
 	\end{aligned} 
 \end{equation} 
 Let 
 \[
\mathrm{I}_1 \coloneqq \big\{ i\in\mathrm{I}\st |\tlam_i|\le \gamma h \big\} \, ,  \quad \mathrm{I}_2\coloneqq \mathrm{I}\setminus\mathrm{I}_1.
 \]
 Clearly, $\mathrm{I}_1\cap\mathrm{I}_2$ is empty and $\mathrm{I}=\mathrm{I}_1\cup\mathrm{I}_2$. 
 Moreover,  $(x\m,\lambda,\xi)\in \Omega$ implies $|\tlam|\le \kappa_1$; as a consequence, for every $i\in\mathrm{I}_2$, $-|\tlam_i|\le - \frac{1}{\kappa_1}|\tlam_i|^2$. Then, by using $\gamma<\bar\gamma_{20}$,  we get
 \begin{align*}
 	&\frac{1}{8}\gamma^2  M k_2^2 |\tlam_{\mathrm{I}}|^2 -\frac{1}{2} \sum_{i \in {\mathrm{I}}} \min \left\{|\tlam_i|^2, \gamma h|\tlam_i|\right\}
 	\\
 	&\qquad \le \sum_{i\in\mathrm{I}_1} \left( \frac{1}{8}\gamma^2 Mk_2^2 - \frac{1}{2} \right) |\tlam_i|^2 + \sum_{i\in\mathrm{I}_2} \left(\frac{1}{8}\gamma^2 Mk_2^2 \kappa_1 -\frac{1}{2}\gamma h \right) |\tlam_i| 
 	\\
 	&\qquad \le -\frac{1}{4}\sum_{i\in\mathrm{I}_1}|\tlam_i|^2 -\frac{1}{4}\gamma h\sum_{i\in\mathrm{I}_2}|\tlam_i| \le -\frac{1}{4}\sum_{i\in\mathrm{I}_1}|\tlam_i|^2 -\frac{1}{4}\gamma \frac{h}{\kappa_1}\sum_{i\in\mathrm{I}_2}|\tlam_i|^2\\
 	&\qquad \le -\frac{1}{4}\min\left\{1,\, \gamma \frac{h}{\kappa_1}\right\}|\tlam_{\mathrm{I}}|^2 \, .
 \end{align*}
Thus, from \eqref{e.V_le_eps} and by using \eqref{e.bounds_P}, \eqref{d.V}, and Item~\ref{item.lem.Mp.sandwitch} of Lemma~\ref{lem.MP}, we obtain
\begin{equation}\label{e.V_le_epsilon}
	\begin{aligned}
        V(x\m^+,\lambda^+,\xi^+)  	&\le  V(x\m,\lambda,\xi)   - \frac{3}{16 c_u^2} V_{\rm net}(\xi)\\
        &\qquad-  \frac{1}{12}\min\left\{  \gamma \mu_f,\,\gamma^2{M} k_2^2,\, 2 \gamma\frac{h}{\kappa_1},\, 2     \right\} V_{\rm opt}(x\m,\lambda) 
        \\
        &\le V(x\m,\lambda,\xi) -  \omega_\gamma V(x\m,\lambda,\xi)
	\end{aligned}
\end{equation}
where 
\begin{equation}\label{d.omega_gamma}
	\omega_\gamma \coloneqq  \min\left\{ \frac{1}{12}\gamma \mu_f, \, \frac{1}{12}\gamma^2 k_2^2,\, \frac{1}{6}\gamma\frac{h}{\kappa_1},\, \frac{1}{6},\, \frac{1}{8c_u}  \right\} \in (0,1),
\end{equation}
and the last inequality follows from the fact that $M \ge 1$ according to \eqref{d.M_beta}.

The following lemma summarizes the results obtained in the present section regarding the increment of $V$ along the solutions of the core subsystem. 
\begin{lemma}\label{lem.V-Omega}
	Suppose that Assumptions~\ref{ass.optimization_problem_data}, \ref{ass.regular_θ*}, and \ref{ass.network} hold.
	Then,  
	\begin{equation}\label{e.V_in_Omega}
		\forall  (x\m,\lambda,\xi)\in \Omega,\quad V(x\m^+ ,\lambda^+ ,\xi^+ )\le V(x\m,\lambda,\xi) - \min\left\{ \frac{1}{8}\gamma\mu_f\varepsilon^2,\, \omega_\gamma V(x\m,\lambda,\xi) \right\}.
	\end{equation}
	In particular, $\Omega$ is forward-invariant for the core subsystem~\eqref{s.reduced}. 
\end{lemma}
\begin{proof}
	Inequality \eqref{e.V_in_Omega} directly follows by joining \eqref{e.V_g_epsilon} for $|x\m-\theta\sr|>\varepsilon$ and \eqref{e.V_le_epsilon} for $|x\m-\theta\sr|\le\varepsilon$.
	Invariance of $\Omega$ is a trivial consequence of \eqref{e.V_in_Omega}.
    \end{proof} 
    
Lemma~\ref{lem.V-Omega} implies the sought stability properties of the core subsystem. Indeed, by \eqref{d.kappa0-eq} and  Lemma~\ref{lem.Omega}, for every initial condition $(x^0,z^0,\lambda^0)\in\Xi_0$, we have that $(x\m^0,\lambda^0,\xi^0)\in \Omega$. Since Lemma~\ref{lem.V-Omega} implies forward-invariance of $\Omega$, it therefore also implies uniform boundedness of the trajectories of Algorithm~\eqref{s.algorithm} originating from $\Xi_0$.
Moreover, it also implies convergence of $(x,z,\lambda)$ to the optimal steady state $(\1\theta\sr,z\sr,\lambda\sr)$, since $V$ is a decreasing Lyapunov function. The next section refines the result by providing the sought exponential bound, thereby concluding the proof.

\subsection{Convergence Rate and Exponential Bound} \label{sec.proof.convergence_rate}

This section complements such results by proving the exponential bound \eqref{e.exp-bound}, thereby also providing an estimate of the convergence rate. This section follows similar arguments of~\cite[Sec. 4.7]{bin2024semiglobal}. 

First, the following lemma shows that every solution of Algorithm~\ref{s.algorithm} enters uniformly in each neighborhood of the optimal point.
\begin{lemma}\label{lem.unif-conv}
	Suppose that Assumptions~\ref{ass.optimization_problem_data}, \ref{ass.regular_θ*}, and \ref{ass.network} hold.
	Moreover, for each $\eta>0$, let
	\begin{equation}\label{d.T-eta}
		T_\gamma(\eta) \coloneqq \frac{3\kappa_0^2 - \min\{\varepsilon^2,\eta^2\}}{\omega_\gamma\min\{\varepsilon^2,\eta^2\}} .
	\end{equation}
	Then, every solution of \eqref{s.algorithm} originating from $\Xi_0$ satisfies
	\begin{equation*}
		\forall t\ge T_\gamma(\eta),\qquad V(x\m^t,\lambda^t,\xi^t)\le \frac{1}{2}\min\{\eta^2,\varepsilon^2\}.
	\end{equation*}
\end{lemma}
\begin{proof}
    See Appendix~\ref{proof.lem.unif-conv}.
\end{proof}

Now, pick a solution $(x,z,\lambda)$ of \eqref{s.algorithm} originating from $\Xi_0$.
In view of Lemma~\ref{lem.V-Omega}, Lemma~\ref{lem.unif-conv} applied with 
\begin{equation}\label{d.proof.eta2}
	\eta^2\coloneqq \frac{\gamma\mu_f \varepsilon^2}{4\omega_\gamma} 
\end{equation}
implies that
\begin{equation*}
	\forall t\ge T_\gamma(\eta), \quad V(x\m^{t+1},\lambda^{t+1},\xi^{t+1})\le (1-\omega_\gamma) V(x\m^t,\lambda^t,\xi^t)\, .
\end{equation*}
Since $1-\omega_\gamma>0$ and $t\mapsto V(x\m^t,\lambda^t,\xi^t)$ is decreasing in view of Lemma~\ref{lem.V-Omega}, we get:
\begin{equation}\label{e.exp-V-1}
	\begin{aligned}
		\forall t\ge  T_\gamma(\eta), \quad V(x\m^t,\lambda^t,\xi^t)&\le   (1-\omega_\gamma)^{t-T_\gamma(\eta)}  V(x\m^{T_\gamma(\eta)},\lambda^{T_\gamma(\eta)},\xi^{T_\gamma(\eta)})
		\\
		&\le (1-\omega_\gamma)^{t-T_\gamma(\eta)}V(x\m^0,\lambda^0,\xi^0)\, .
	\end{aligned} 
\end{equation}
Since $1-\omega_\gamma\in(0,1)$, for all $t<T_\gamma(\eta)$ one has $(1-\omega_\gamma)^{t-T_\gamma(\eta)}>1$. Hence, by Lemma~\ref{lem.V-Omega}, it follows that 
\begin{equation}\label{e.exp-V-2}
	\forall t\in\{0,\ldots,T_\gamma(\eta)-1\},\quad V(x\m^t,\lambda^t,\xi^t)\le V(x\m^0,\lambda^0,\xi^0) \le (1-\omega_\gamma)^{t-T_\gamma(\eta)} V(x\m^0,\lambda^0,\xi^0).
\end{equation}
Joining \eqref{e.exp-V-1} and \eqref{e.exp-V-2},  finally yields
\begin{equation} \label{e.exp-V-3}
	\forall t\in\N, \quad V(x\m^t,\lambda^t,\xi^t)\le   (1-\omega_\gamma)^{t-T_\gamma(\eta)}  V(x\m^0,\lambda^0,\xi^0) \, .
\end{equation}
Then, the claim of the theorem follows with
\begin{align}\label{d.proof.c_mu}
	c&\coloneqq \sqrt{n \max\left\{ 3,\, \frac{c_u}{c_l}\right\} (1-\omega_\gamma)^{-T_\gamma(\eta)}}, & \mu &\coloneqq \sqrt{1-\omega_\gamma} ,
\end{align}
 directly  from \eqref{e.sandwitch_dist_A},  \eqref{e.exp-V-3} and by the fact that, in view of Lemma~\ref{lem.MP} and \eqref{e.bounds_P}, for all $(x\m,\lambda,\xi)\in\Omega$, the following inequality hods:
\begin{align*}
\frac{1}{2}\min\left\{1, 3\frac{c_l}{c_u}\right\} |(\tx\m,\tlam,\txi)|^2	\le V(x\m,\lambda,\xi) \le \frac{3}{2} |(\tx\m,\tlam,\txi)|^2.
\end{align*}

 \section{Conclusions}
 This paper dealt with a distributed primal-dual algorithm for constrained consensus optimization over a network of agents.
 The algorithm fuses the centralized primal-dual gradient method analyzed in \cite{bin2024semiglobal} and the distributed unconstrained method of \cite{wang_control_2010,wang2011control,binStabilityLinearConvergence2022}.
 The analysis rationale is based on the  separation between the optimization and the network dynamics; the former refers to the evolution of the average estimate of the agents toward the optimal solution, the latter refers to the evolution of the consensus error over the network.
 It is found that, indeed, the considered algorithm can be naturally decomposed in the mutual interconnection of these dynamics, which take place at a different time scale. 
 These two dynamics are mutually coupled, and this makes the stability analysis rather challenging, since it requires suitable time-scale separation, Lyapunov analysis, and small-gain tools from control theory.
 
 Remarkably, the analysis of the optimization dynamics leveraged the result of~\cite{bin2024semiglobal}, which was taken off the shelf and generalized to the present context.
 Indeed, a relevant ontological consequence of the analysis approach proposed here is that the network effect, which is what makes the algorithm studied here different from its centralized version \cite{bin2024semiglobal}, can be seen as a fast parasite dynamics that, for small enough $\gamma$, can be simply ignored as robustness of semiglobal exponential stability with respect to such kinds of  unmodeled dynamics is already implied by the nominal result.
 In turn, this rationale enables the possibility to handle many of the imaginable extensions of the algorithm, for instance to a class of time-varying networks with dwell-time constraints or small-delayed communication, as these can be seen as a perturbed version of the centralized algorithm of \cite{bin2024semiglobal} and they can therefore be  accounted for by the robustness of semiglobal exponential stability provided by 
\cite{bin2024semiglobal} for the centralized algorithm.

\appendix
\section{Proof of Lemma \ref{lem.equilibria}}\label{proof.lem.equilibria} 
	Let $(x\eq,z\eq,\lambda\eq)$ be an equilibrium of~\eqref{s.algorithm}.
	Equation \eqref{s.algorithm.z} implies $Kx\eq=0$ that, in view of \eqref{d.K}, implies $x\eq = \1 \theta\sr$ for some $\theta\sr\in\R$.  
	Equation \eqref{s.algorithm.x} implies $Kz\eq = - \gamma{n} \nablaL(\1\theta\sr,\lambda\eq)$. Pre-multiplying both sides by $\1\T$ and using \eqref{d.K} yields $0 = \sum_{i\in\cN}\nabla\ell_i(\theta\sr,\lambda\eq_i)$, which is \eqref{d.KKT.gradients} with $\lambda\sr=\lambda\eq$.
	Equation \eqref{d.KKT.constraints}  with $\lambda\sr=\lambda\eq$  directly follows from \eqref{s.algorithm.λ}. 
	Finally, from \eqref{s.algorithm.λ} we obtain $0 = \max\{-\lambda\eq_{i,j},\, \gamma g_{i,j}(\theta\sr)\}$ for all $i\in\cN$, $j=1,\dots, m_i$.
	Therefore,  if $\lambda\eq_{i,j}>0$, then necessarily $g_{i,j}(\theta\sr)=0$ and \eqref{d.KKT} with $\lambda\sr=\lambda\eq$ holds true.  Conversely, if $\lambda\eq_{i,j}=0$, \eqref{d.KKT.active_constraints} trivially holds as well. Hence, $(x\eq,\lambda\eq)=(\1\theta\sr,\lambda\sr)$ with $(\theta\sr,\lambda\sr)$ solves \eqref{d.KKT}.
	
	As for the second claim, let $(\theta\sr,\lambda\sr)$ satisfy \eqref{d.KKT} and define $z\eq\coloneqq  -\gamma {n} {S} (S\T K S)\inv S\T \Phi(\1\theta\sr,\lambda\sr)$, $x\eq\coloneqq \1\theta\sr$, and $\lambda\eq\coloneqq \lambda\sr$.
	As $x\eq\in\linspan\1 = \ker K$ (in view of \eqref{d.K}), then the constant trajectory $t\mapsto(z\eq,x\eq)$ satisfies \eqref{s.algorithm.z}.
	Moreover,  
    \eqref{d.KKT.constraints}-\eqref{d.KKT.active_constraints} imply $0=\max\{-\lambda\sr_{i,j}, \, \gamma g_{i,j}(\theta\sr)\}$ for all $i\in\cN$, $j=1,\dots, m_i$, from which it follows that $t\mapsto (\lambda\eq,x\eq)$ satisfies \eqref{s.algorithm.λ}.
	Finally, in view of \eqref{d.K}, we have 
	\begin{align*}
		&Kz\eq = {SS\T K}z\eq = -\gamma {n}SS\T K S (S\T K S)\inv S\T \Phi(\1\theta\sr,\lambda\sr) \\&= - \gamma {n}S S\T \Phi(\1\theta\sr,\lambda\sr) =-\gamma {n} \Phi(\1\theta\sr,\lambda\sr),
	\end{align*}
	where the last equality follows from the fact that, by \eqref{d.KKT.gradients}, $ \Phi(\1\theta\sr,\lambda\sr)\in\ker(\1\T)$. Hence, $t\mapsto(x\eq,z\eq,\lambda\eq)$ satisfies \eqref{s.algorithm.x}.

\section{Proof of Lemma~\ref{lem.Omega}}\label{proof.lem.Omega}
Regarding Item \ref{item.lem.Omega.sandwich_Omega},
pick $(x\m,\lambda,\xi)\in \overline{\mathbb{B}}_{\kappa_0}(\theta\sr,\lambda\sr,\xi\sr)$. Then, $(x\m,\lambda,\xi)\in\cK$ and  $(x\m,\lambda)\in \overline{\mathbb{B}}_{\kappa}(\theta\sr,\lambda\sr)$, where, we recall, $\kappa$ has been fixed in Section~\ref{sec.proof.preliminary} to satisfy~\eqref{d.kappa}, and $\cK$ is defined in~\eqref{set_k}.
Hence, it follows from Item~\ref{item.lem.Mp.sandwitch} of Lemma~\ref{lem.MP} and from \eqref{e.bounds_P}  that
\begin{align*}
	V(x\m,\lambda,\xi)
 &=V_{\rm opt}(x\m, \lambda)+ {\frac{3}{2c_u}}V_{\rm net}(\xi)
 \\&\le \frac{3}{2}|(x\m-\theta\sr,\lambda-\lambda\sr)|^2 +  \frac{3}{2} |\xi-\xi\sr|^2 \\&=  \frac{3}{2}  (x\m-\theta\sr,\lambda-\lambda\sr,\xi-\xi\sr)|^2\le  \frac{3}{2} \kappa_0^2, 
\end{align*}
which shows that $\overline{\mathbb{B}}_{\kappa_0}(\theta\sr,\lambda\sr,\xi\sr)\subset\Omega$.

Now, pick $(x\m,\lambda,\xi)\in \Omega$. From   Item~\ref{item.lem.Mp.sandwitch} of Lemma~\ref{lem.MP} (notice that $(x\m,\lambda,\xi)\in \Omega$ implies $(x\m,\lambda)\in\overline{\mathbb{B}}_{\kappa}(\theta\sr,\lambda\sr)$) and  \eqref{e.bounds_P}, 
we obtain
\begin{align*}
	V(x\m,\lambda,\xi) &\ge \frac{1}{2}|(x\m-\theta\sr,\lambda-\lambda\sr)|^2 + \frac{3c_l}{2c_u} |\xi-\xi\sr|^2  
	\\
	&\ge   \min\left\{\frac{1}{2} ,\, \frac{3c_l}{2c_u} \right\}|(x\m-\theta\sr,\lambda-\lambda\sr,\xi-\xi\sr)|^2.
\end{align*} 
Thus,
\begin{align*}
		|(x\m-\theta\sr,\lambda-\lambda\sr,\xi-\xi\sr)|^2\le    \max\left\{2,\ \frac{2c_u}{3c_l} \right\}	V(x\m,\lambda,\xi) \le   \max\left\{3,\ \frac{c_u}{c_l} \right\} \kappa_0^2 = \kappa_1^2 ,
\end{align*}
where $\kappa_1$ is defined in \eqref{d.kappa1},
which proves $\Omega\subset \overline{\mathbb{B}}_{\kappa_1}(\theta\sr,\lambda\sr,\xi\sr)$.

Regarding Item \ref{item.lem.Omega.plus},  pick $(x\m,\lambda,\xi)\in\overline{\mathbb{B}}_{\kappa_1}(\theta\sr,\lambda\sr,\xi\sr)$ and, with reference to the core subsystem \eqref{s.reduced}, consider the notation~\eqref{e.notation_plus}. By using \eqref{d.k5_k6_k7} and $\gamma<\bar\gamma_{13}$, proceeding as in \cite[Eq. (22)]{bin2024semiglobal} 
we obtain
\begin{align*} 
		\left|\tx\m^+\right|^2 
		& =|\tilde{x}\m|^2-2 \gamma  \tx\m    1\T  \nablaL(x, \lambda) +\gamma^2| 1\T \nablaL(x, \lambda)|^2  \\&  \leq |\tilde{x}\m|^2+ {2}\gamma \kappa_1 k_5+ {\gamma^2} k_5^2
		\\& \leq |\tilde{x}\m|^2+\frac{1}{3}  ,
		\\
		|\tlam^+|^2 
		&\leq |\tlam|^2 +2\gamma\sprod{\tlam}{g(x)} + \gamma^2 |g(x)|^2  
	\\& 	\leq|\tlam|^2 + 2\gamma \kappa_1 k_7  +\gamma^2 k_7^2
		\\& \leq|\tilde{\lambda}|^2+\frac{1}{3} .
\end{align*} 
 Besides, $\gamma<\bar\gamma_{14}$ and \eqref{d.k1_k2_k3_k4}   imply 
\begin{equation}\label{e.proof.Omega.Vnet}
	\begin{aligned}
		V_{\rm net}(\xi^+) & \leq V_{\rm net}(\xi) -\frac{1}{2}|\txi|^2 +  \gamma^2 n^2 \big(2|A\T PB|^2+|B\T PB|\big)  |\Phi(\1 \theta\sr,\lambda\sr)-\Phi(x,\lambda)|^2 
		\\
		&\le 
    V_{\rm net}(\xi) -\frac{1}{2}|\txi|^2 +  \gamma^2  n^2 k_0^2 \big(2|A\T PB|^2+|B\T PB|\big)  \big( |\tx\m|+|\txi| + |\tlam|\big)^2 
			\\
		&\le 
    V_{\rm net}(\xi) + \left(\gamma^2 n^2 k_0^2  \big(6|A\T PB|^2+3|B\T PB|\big) -\frac{1}{2}\right)|\txi|^2 \\
		&\qquad\quad+   3\gamma^2 n^2 k_0^2 \big(2|A\T PB|^2+|B\T PB|\big) \big( |\tx\m|^2  + |\tlam|^2\big) 
		\\
		&\le V_{\rm net}(\xi) -\frac{1}{4}|\txi|^2+   3\gamma^2 n^2 k_0^2 \big(2|A\T PB|^2+|B\T PB|\big) \big( |\tx\m|^2  + |\tlam|^2\big),
	\end{aligned}
\end{equation}
where the first inequality is obtained from \eqref{e.V_net_1} by substituting $v$ with \eqref{d.v}.
By combining $|\tx\m|^2  + |\tlam|^2\le \kappa_1^2{< \kappa^2}$ and   $\gamma<\bar\gamma_{14}$, we get
\begin{align*}
		V_{\rm net}(\xi^+)   \leq V_{\rm net}(\xi)+  \frac{c_l}{3} . 
\end{align*} 
Hence, from \eqref{e.bounds_P} it follows that
\begin{align*}
	|\txi^+|^2 &\leq \frac{1}{c_l} V_{\rm net}(\xi^+)
	\leq 
	\frac{1}{c_l} \left( V_{\rm net}(\xi) + \frac{c_l}{3}\right)
	 \leq \frac{c_u}{c_l} |\txi|^2 +  \frac{1}{3} \, .
\end{align*}
Owing to the previous bounds for $|\tx\m^+|^2$, $|\tlam^+|^2$, and $|\txi^+|^2$, we obtain
 \begin{equation*}
	\begin{aligned}
		|\tx\m^+|^2 +  |\tlam^+|^2  + |\txi^+|^2  &\leq  \frac{c_u}{c_l} (  |\tx\m|^2 +  |\tlam|^2  + |\txi |^2  + 1 ) \\
		&\le \frac{c_u}{c_l}(\kappa_1^2+1) \le \left(\sqrt{\frac{c_u}{c_l}}(\kappa_1+1)\right)^2 \le \kappa^2 \, ,
	\end{aligned}
\end{equation*}
which implies 
$( x\m^+, \lambda^+, \xi^+  ) \in \mathcal{K}$.

Finally, Item~\ref{item.lem.Omega.opt}, follows by noticing that $V_{\rm opt}(x\m,\lambda)\le V(x\m,\lambda,\xi)$ for every $(x\m,\lambda,\xi)\in\R\x(\R_{\ge 0})^m \x \R^{2(n-1)}$.

\section{Proof of Lemma \ref{lem.W}}\label{proof.lem.W}
First, we expand $W_\gamma$ as follows
\begin{align*}
	W_\gamma&(x\m,\lambda,w,d)=    2\gamma \left(\tx\m-\gamma\1\T\Phi(\1 x\m,\lambda)\right) w + \gamma^2 w^2 
	\\
	& \qquad  + |\max\{-\lambda\sr,\tlam  + \gamma g(\1 x\m)+\gamma d \}|^2  
	- |\max\{-\lambda\sr,\tlam  + \gamma g(\1 x\m)\}|^2
	+ \cI \, , 
\end{align*}
where $\cI$ is defined as
\begin{align*}
	\cI &\coloneqq  \gamma \beta \left( \tx\m-    \gamma\1\T\Phi(\1 x\m,\lambda)+\gamma w\right)  \nabla g\left(\1(x\m -  \gamma\1\T\Phi(\1 x\m,\lambda) +\gamma w ) \right)
	\\& \qquad\qquad  \max\{-\lambda\sr,\tlam  + \gamma g(\1x\m) +\gamma d \}
	\\
	&\quad - \gamma \beta \left( \tx\m-   \gamma\1\T\Phi(\1 x\m,\lambda) \right) \nabla g\left(\1(x\m -
	\gamma\1\T\Phi(\1 x\m,\lambda) )\right)
	\\&\qquad\qquad  \max\{-\lambda\sr,\tlam  + \gamma g(\1x\m)\}   \, .
\end{align*}
{For every $a,b,c\in\R$, the following facts hold }
\begin{align*}
    & a^2-b^2 =   (a-b)^2 +2b(a-b)
    \\
    &|\max\{a,b\}-\max\{a,c\}|\le |b-c|.
\end{align*}
Hence, we can write
\begin{align*}
	&|\max\{-\lambda\sr,\tlam  + \gamma g(\1 x\m)+\gamma d \}|^2  
	- |\max\{-\lambda\sr,\tlam  + \gamma g(\1 x\m)\}|^2 
	\\
	&\  = \left(|\max\{-\lambda\sr,\tlam  + \gamma g(\1 x\m)+\gamma d \}|- |\max\{-\lambda\sr,\tlam  + \gamma g(\1 x\m) \} |\right)^2 
    \\
	&\  \quad 
    + 2 |\max\{-\lambda\sr,\tlam  + \gamma g(\1 x\m) \} |\times \\
    &\ \quad \times \left(|\max\{-\lambda\sr,\tlam  + \gamma g(\1 x\m)+\gamma d \}|- |\max\{-\lambda\sr,\tlam  + \gamma g(\1 x\m) \}| \right)\\
	&\le \gamma^2 |d|^2 +2\gamma   |d| \big|  \max\{-\lambda\sr,\tlam  + \gamma g(\1 x\m) \} \big|.
\end{align*}
Moreover, for every $a,b,\bar a,\bar b\in\R$, we have
\begin{align*}
    ab-\bar a \bar b = (a-\bar a)(b-\bar b)+\bar b(a-\bar a)+\bar a(b-\bar b).
\end{align*}
By using this equality with
$a= ( \tx\m-    \gamma\1\T\Phi(\1 x\m,\lambda)+\gamma w)   \nabla g(\1(x\m -  \gamma\1\T\Phi(\1 x\m,\lambda) +\gamma w ) )$, 
$b=\max\{-\lambda\sr,\tlam  + \gamma g(\1x\m) +\gamma d \}$, 
$\bar a=( \tx\m- \gamma\1\T\Phi(\1 x\m,\lambda) )  \nabla g(\1(x\m -
\gamma\1\T\Phi(\1 x\m,\lambda) ))$, 
and $\bar b= \max\{-\lambda\sr,\tlam  + \gamma g(\1x\m)\}$, 
we can rewrite $\cI$ as
\begin{align*}
	\cI &= \gamma\beta  \Big[  \left( \tx\m-    \gamma\1\T\Phi(\1 x\m,\lambda)+\gamma w\right)   \nabla g\left(\1(x\m -  \gamma\1\T\Phi(\1 x\m,\lambda) +\gamma w ) \right) \\
	&\qquad\qquad\qquad - \left( \tx\m-   \gamma\1\T\Phi(\1 x\m,\lambda) \right)  \nabla g\left(\1(x\m -
	\gamma\1\T\Phi(\1 x\m,\lambda) )\right)\Big] \cdot 
 \\ &\hspace{1.5cm}\cdot\Big[  \max\{-\lambda\sr,\tlam  + \gamma g(\1x\m) +\gamma d \} - \max\{-\lambda\sr,\tlam  + \gamma g(\1x\m)\}\Big]  \\
	&\quad + \gamma\beta     \max\{-\lambda\sr,\tlam  + \gamma g(\1x\m)\} \cdot  \\
	&\qquad\qquad  \cdot \Big[  \left( \tx\m-    \gamma\1\T\Phi(\1 x\m,\lambda)+\gamma w\right)   \nabla g\left(\1(x\m -  \gamma\1\T\Phi(\1 x\m,\lambda) +\gamma w ) \right) \\
	&\qquad\qquad\qquad - \left( \tx\m-   \gamma\1\T\Phi(\1 x\m,\lambda) \right)  \nabla g\left(\1(x\m -
	\gamma\1\T\Phi(\1 x\m,\lambda) )\right)\Big]   \\
	&\quad + \gamma\beta \left( \tx\m-   \gamma\1\T\Phi(\1 x\m,\lambda) \right)  \nabla g\left(\1(x\m -
	\gamma\1\T\Phi(\1 x\m,\lambda) )\right) \cdot \\
	&\qquad\qquad  \cdot \Big[  \max\{-\lambda\sr,\tlam  + \gamma g(\1x\m) +\gamma d \} - \max\{-\lambda\sr,\tlam  + \gamma g(\1x\m)\}\Big]
	\\
	&\le \gamma^2\beta |d|   |\cJ_1| + \gamma\beta|\cJ_1| |\max\{-\lambda\sr,\tlam+\gamma g(\1\tx\m)\}|    + \gamma^2 \beta |d|  |\cJ_2|,
\end{align*}
where, in the last inequality, we used again the previous algebraic bounds, the \linebreak Cauchy–Schwarz inequality, and we defined
\begin{align*}
	\cJ_1 &\coloneqq  \left( \tx\m-    \gamma\1\T\Phi(\1 x\m,\lambda)+\gamma w\right)   \nabla g\left(\1(x\m -  \gamma\1\T\Phi(\1 x\m,\lambda) +\gamma w ) \right) \\
	&\qquad\qquad\qquad - \left( \tx\m-   \gamma\1\T\Phi(\1 x\m,\lambda) \right)  \nabla g\left(\1(x\m -
	\gamma\1\T\Phi(\1 x\m,\lambda) )\right) \, ,
 \\
	\cJ_2&:=\left( \tx\m-   \gamma\1\T\Phi(\1 x\m,\lambda) \right)  \nabla g\left(\1(x\m -
	\gamma\1\T\Phi(\1 x\m,\lambda) )\right).
\end{align*}
Thus, for every $(x\m,\lambda,w,d)\in\overline{\mathbb{B}}_{\kappa_1}(\theta\sr,\lambda\sr)
\x\R \x \R^m$, we can write
\begin{equation}\label{e.proof.W.W_1}
	\begin{aligned}
		W_\gamma (x\m,\lambda,w,d)&\le  2\gamma \left(\tx\m-\gamma\1\T\Phi(\1 x\m,\lambda)\right) w + \gamma^2 \big(|w|^2+|d|^2 + \beta |d|   (|\cJ_1|+|\cJ_2|) \big)
		\\
		&\qquad 
		+ \gamma  |  \max\{-\lambda\sr,\tlam  + \gamma g(\1 x\m) \} |  \big(  2 |d|  +  \beta  |\cJ_1| \big)  .
	\end{aligned}
\end{equation}
{Now, let us specialize the previous bound for the specific choice of $w$ and $d$ given in~\eqref{d.w_d}, i.e., set}
\begin{align}\label{d.proof.W.w_d}
	w  &=   \1\T \big( \nablaL(\1x\m ,\lambda )-\nablaL(x,\lambda)\big) ,&
	d &=  g(x)-g(\1x\m).
\end{align}
Then, since $\overline{\mathbb{B}}_{\kappa_1}(\theta\sr,\lambda\sr,\xi\sr)\subset\cK$, by using \eqref{d.k1_k2_k3_k4}, we obtain 
\begin{align}\label{e.proof.W.wd}
	  |w|&\le |\1 | k_0 |x\p| \le \sqrt{n} k_0 |\txi|,&
	  |d|&\le k_1|x\p| \le k_1 |\txi|, 
\end{align}
for every $(x\m,\lambda,\xi)\in\overline{\mathbb{B}}_{\kappa_1}(\theta\sr,\lambda\sr,\xi\sr)$.
As a consequence, by using \eqref{d.k1_k2_k3_k4} again, for every $(x\m,\lambda,\xi)\in\overline{\mathbb{B}}_{\kappa_1}(\theta\sr,\lambda\sr,\xi\sr)$, {the following holds}
\begin{equation}\label{e.proof.W.term1}
	2\gamma \left(\tx\m-\gamma\1\T\Phi(\1 x\m,\lambda)\right) w \le 2\gamma \sqrt{n}k_0 |\tx\m||\txi| + 2\gamma^2 \sqrt{n}k_0k_2 (|\tx\m|+|\tlam|)|\txi|. 
\end{equation}
Moreover, for every $i\in\cN$ and $j=1,\dots,m_i$, we can write
\begin{equation*}
 \max\big\{ -\lambda_{i,j}\sr,\,\tlam_{i,j}  + \gamma g_{i,j}(x\m) \big\} = 
  \begin{cases}
		\lambda_{i,j} + \gamma g_{i,j}(x\m) - \lambda\sr_{i,j} \, , & \text{if } \lambda_{i,j} + \gamma g_{i,j}(x\m)  \geq 0  \, , \\
		- \lambda\sr_{i,j} \, ,  & \text{if }  \lambda_{i,j} + \gamma g_{i,j}(x\m)  < 0 \, .
	\end{cases} 
\end{equation*}
Since $\lambda_{i,j} + \gamma g_{i,j}(x\m)  < 0$ implies $\lambda\sr_{i,j} < -(\tlam_{i,j} + \gamma g_{i,j}(x\m))\le |\tlam_{i,j} + \gamma g_{i,j}(x\m)|$, and since  $\lambda\sr_{i,j}\ge 0$, we then deduce
\begin{equation*}
	|\max\big\{ -\lambda_{i,j}\sr,\, \tlam_{i,j}  + \gamma g_{i,j}(x\m) \big\}| \le |\tlam_{i,j} + \gamma g_{i,j}(x\m)|.
\end{equation*}
On the one hand, for every $(i,j)\in \mathrm{\Lambda}$ (i.e., $(i,j)$ is associated with an active constraint),  $g_{i,j}(\theta\sr) = 0$. Thus, by using \eqref{d.k1_k2_k3_k4} and the previous inequality, for every $(x\m,\lambda,\xi)\in\overline{\mathbb{B}}_{\kappa_1}(\theta\sr,\lambda\sr,\xi\sr)$, one gets
\begin{equation*}
	 |\max\{-\lambda\sr_{i,j},\,\tlam_{i,j}  + \gamma g_{i,j}(x\m)\}|  \leq |  \tlam_{i,j} | + \gamma | g_{i,j}(x\m) - g_{i,j}(\theta\sr)  | \le  |  \tlam_{i,j} | +   \gamma k_1 |\tx\m |.
\end{equation*} 
On the other hand, for every $(i,j)\in \mathrm{I}$ (i.e., $(i,j)$ is associated with an inactive constraint),  $g_{i,j}(\theta\sr) < 0$ and $\lambda\sr_{i,j}=0$. Hence, $\tlam_{i,j}=\lambda_{i,j}\ge 0$, and $-g_{i,j}(\theta\sr)>0$.
If $g_{i,j}(x\m)\le 0$, then
\begin{align*}
	\big|\max\big\{-\lambda\sr_{i,j},\,\tlam_{i,j}  + \gamma g_{i,j}(x\m)\big\}\big| 
	&=\begin{cases}
		\lambda_{i,j}+\gamma g_{i,j}(x\m) \, , & 	\text{if}\ \lambda_{i,j}+\gamma g_{i,j}(x\m) \ge 0 \, ,\\
		0  \, ,& \text{if}\	\lambda_{i,j}+\gamma g_{i,j}(x\m) <0 \, ,
	\end{cases}\\
	 &\le \lambda_{i,j} = |\lambda_{i,j}|=|\tlam_{i,j}| \, .
\end{align*}
If, instead, $g_{i,j}(x\m)>0$, we get
\begin{align*}
	\big|\max\big\{-\lambda\sr_{i,j},\,\tlam_{i,j}  + \gamma g_{i,j}(x\m)\big\}\big| 
	&= \lambda_{i,j}  + \gamma g_{i,j}(x\m) 
	\le \lambda_{i,j}  + \gamma g_{i,j}(x\m) - \gamma g_{i,j}(\theta\sr)\\
	&\le |\tlam_{i,j}| + \gamma k_1 |\tx\m| ,
\end{align*}
for   every $(x\m,\lambda,\xi)\in\overline{\mathbb{B}}_{\kappa_1}(\theta\sr,\lambda\sr,\xi\sr)$.
Therefore, for all $i\in\cN$, $j=1,\dots,m_i$, and every 
$(x\m,\lambda,\xi)\in\overline{\mathbb{B}}_{\kappa_1}(\theta\sr,\lambda\sr,\xi\sr)$,  one has
\begin{equation*}
		\big|\max\big\{-\lambda\sr_{i,j},\,\tlam_{i,j}  + \gamma g_{i,j}(x\m)\big\}\big| \le |\tlam_{i,j}| + \gamma k_1 |\tx\m| 
\end{equation*}
and, hence,
\begin{equation}\label{e.proof.W.max}
	\begin{aligned}
		\big|\max\big\{-\lambda\sr,\,\tlam + \gamma g(\1 x\m)\big\}\big| &\le \sum_{(i,j)\in\mathrm{\Lambda}\cup\mathrm{I}} \big|\max\big\{-\lambda\sr_{i,j},\,\tlam_{i,j}  + \gamma g_{i,j}(x\m)\big\}\big|
		\\
		&\le m(|\tlam| + \gamma k_1 |\tx\m|) \, ,
	\end{aligned}
\end{equation}
where we recall that $m$ is the total number of constraints.
	
Next, by substituting \eqref{d.proof.W.w_d} into the expression of $\cJ_1$, we obtain 
\begin{align*}
	\cJ_1 & =  \left( \tx\m-    \gamma\1\T\Phi(x,\lambda)\right)   \nabla g\left(\1(x\m -  \gamma\1\T\Phi(x,\lambda))\right) \\
	&\qquad\qquad\qquad - \left( \tx\m-   \gamma\1\T\Phi(\1 x\m,\lambda) \right)  \nabla g\left(\1(x\m -
	\gamma\1\T\Phi(\1 x\m,\lambda) )\right).
\end{align*}
By Lemma~\ref{lem.Omega}, $(x\m,\lambda,\xi)\in\overline{\mathbb{B}}_{\kappa_1}(\theta\sr,\lambda\sr,\xi\sr)$ implies $(x\m^+,\lambda^+,\xi^+)\in \cK$. Thus, we can use \eqref{d.k1_k2_k3_k4} and \eqref{d.k5_k6_k7} to bound  $\nabla g(\1(x\m-    \gamma\1\T\Phi(x,\lambda))) = \nabla g(\1 x\m^+)$  in the expressions of $\cJ_1$.
Moreover, we can use the same bounds also to bound the term $\nabla g(\1(x\m- \gamma\1\T\Phi(\1 x\m,\lambda)))$  in $\cJ_1$ and $\cJ_2$. 
In fact,  $(x\m',\lambda',\xi')\coloneqq(x\m,\lambda,\xi\sr)$ satisfies $x\m'{}^+ = x\m'- \gamma\1\T\Phi(\1 x\m',\lambda')) = x\m- \gamma\1\T\Phi(\1 x\m,\lambda))$, and $(x\m',\lambda',\xi')\in \overline{\mathbb{B}}_{\kappa_1}(\theta\sr,\lambda\sr,\xi\sr)$ since 
\begin{equation*}
	|x\m'-\theta\sr|^2 + |\lambda'-\lambda\sr|^2 + |\xi'-\xi\sr| = |x\m-\theta\sr|^2 + |\lambda-\lambda\sr|^2 \le \kappa_1^2.
\end{equation*}
Thus, by {Lemma~\ref{lem.Omega}}, it follows that $( x\m- \gamma\1\T\Phi(\1 x\m,\lambda), \lambda'{}^+, \xi'{}^+ )\in\cK$.

Now,
using $ab-\bar a\bar b = a(b-\bar b)+ \bar b(a-\bar a)$ with $a=\tx\m-    \gamma\1\T\Phi(x,\lambda)$, $b=\nabla g (\1(x\m -  \gamma\1\T\Phi(x,\lambda)) )$, $\bar a=\tx\m-   \gamma\1\T\Phi(\1 x\m,\lambda)$, and $\bar b=\nabla g (\1(x\m -
\gamma\1\T\Phi(\1 x\m,\lambda) ) )$,  together with \eqref{d.k1_k2_k3_k4}-\eqref{d.k5_k6_k7}  yields
\begin{equation}\label{e.proof.W.J1}
	\begin{aligned}
			|\cJ_1| &\le \left|\tx\m-    \gamma\1\T\Phi(x,\lambda) \right|     \big| 
				\nabla g \left(\1(x\m -  \gamma\1\T\Phi(x,\lambda)) \right)
				\\&\hspace{4cm}-
				\nabla g\left(\1(x\m -
				\gamma\1\T\Phi(\1 x\m,\lambda) )\right)
			\big|
			\\&\qquad +\left|\nabla g (\1(x\m -
			\gamma\1\T\Phi(\1 x\m,\lambda) ) ) \right|  \left|    \gamma\1\T\Phi(\1 x\m,\lambda) -    \gamma\1\T\Phi(x,\lambda) \right|\\
			&
			\le (\kappa_1+\gamma k_5) k_4   \gamma |\1\1\T (\Phi(x,\lambda) - \Phi(\1 x\m,\lambda)) | \\
			&\qquad + \gamma k_6 |\1| |\Phi(x,\lambda) - \Phi(\1 x\m,\lambda)|\\
			&\le \Big[ \gamma (n\kappa_1 {k_4} k_0 + k_0{k_6}\sqrt{n}) + \gamma^2 n k_0 k_4k_5 \Big] |x\p|\\&  \le \Big[ \gamma (n\kappa_1  {k_4} k_0 + k_0{k_6}\sqrt{n}) + \gamma^2 n k_0 k_4k_5 \Big] |\txi|\\
			&\le \frac{1}{\beta} |\txi|,
	\end{aligned} 
\end{equation} 
for all $(x\m,\lambda,\xi)\in\overline{\mathbb{B}}_{\kappa_1}(\theta\sr,\lambda\sr,\xi\sr)$ where, in the last inequality, we have used $\gamma<\bar\gamma_{15}$.
Similarly, we can use \eqref{d.k1_k2_k3_k4}-\eqref{d.k5_k6_k7} to bound $\cJ_2$ as
\begin{equation}\label{e.proof.W.J2}
	|\cJ_2| \le k_6 (|\tx\m| + \gamma|\1\T\Phi(\1 x\m,\lambda)|) \le k_6(1+\gamma k_2)|\tx\m| + \gamma k_2k_6 |\tlam| \, ,
\end{equation}
for all $(x\m,\lambda,\xi)\in\overline{\mathbb{B}}_{\kappa_1}(\theta\sr,\lambda\sr,\xi\sr)$.

Plugging \eqref{e.proof.W.wd}, \eqref{e.proof.W.term1}, \eqref{e.proof.W.max}, \eqref{e.proof.W.J1}, and \eqref{e.proof.W.J2} into the bound~\eqref{e.proof.W.W_1} for $W_\gamma$, yields the following bound for $\overline{W}_\gamma$ 
\begin{align*}
	&\overline{W}_\gamma(x\m,\lambda,\xi) \\&\le
	  \gamma\Big[2\sqrt{n}k_0 +\gamma\big( 2\sqrt{n}k_0k_2+\beta k_1k_6 +mk_1(2k_1+1) \big) +\gamma^2  \beta k_1k_2k_6     \Big]|\tx\m||\txi|
	  \\
	  &\qquad+ 
	  \gamma\Big[m(2k_1+1) + 2\gamma \sqrt{n}k_0k_2 + \gamma^2 \beta k_1k_2k_6 \Big]|\tlam||\txi|
	  \\
	  &\qquad+  \gamma^2\Big[nk_0^2+k_1^2+  k_1  \Big]|\txi|^2
	  \\
	  &\le 
	  \gamma^2\Big[nk_0^2+k_1^2+ k_1  \Big]|\txi|^2 + 
	  \gamma  6\sqrt{n}k_0    |\tx\m||\txi|
	  +
	  2\gamma m(2k_1+1)   |\tlam||\txi|,
\end{align*}
for all $(x\m,\lambda,\xi)\in\overline{\mathbb{B}}_{\kappa_1}(\theta\sr,\lambda\sr,\xi\sr)$ where, in the last inequality, we used $\gamma<\bar\gamma_{16}$.

Using twice the  Young's inequality \eqref{d.Young}, the first time with $a=3\sqrt{n} k_0|\txi|$, $b=|\tx\m|$, and $\epsilon=\frac{\mu_f c_u}{4(2c_u+3)}$, and the second with $a=\gamma m(2k_1+1)|\tlam|$, $b=|\txi|$, and $\epsilon=\frac{3}{32 c_u}$, leads to
\begin{align*}
	\overline{W}_\gamma(x\m,\lambda,\xi) &\le \gamma \frac{\mu_f c_u}{4(2c_u+3)} |\tx\m|^2 + \frac{32}{3}\gamma^2 c_u m^2(2k_1+1)^2|\tlam|^2\\
	&\qquad + \left[ \frac{3}{32 c_u} + \gamma\frac{36 nk_0^2(2c_u+3)}{\mu_f c_u} + \gamma^2\Big[nk_0^2+k_1^2+k_1  \Big] \right]|\txi|^2,
\end{align*}
for all $(x\m,\lambda,\xi)\in\overline{\mathbb{B}}_{\kappa_1}(\theta\sr,\lambda\sr,\xi\sr)$.
Using $\gamma<\bar\gamma_{17}$ (cf.~\eqref{d.bargamma_2}) finally yields the result.

 \section{Proof of Lemma~\ref{lem.unif-conv}}
\label{proof.lem.unif-conv}
Fix $\eta>0$ arbitrary and let $\upsilon\coloneqq\frac{1}{2}\min\{ \eta^2,\varepsilon^2\}$.
In order to ease the notation, let $T\coloneqq T_\gamma(\eta)$.
Pick a solution $(x,z,\lambda)$ of Algorithm~\eqref{s.algorithm} originating from $\Xi_0$. 
Let $\tau>0$ be such that $V(x\m^t,\lambda^t,\xi^t)>\upsilon$ for all $t< \tau$ and $V(x\m^\tau,\lambda^\tau,\xi^\tau)\le \upsilon$.
Existence of such a $\tau\ge0$ follows from Lemma~\ref{lem.V-Omega}, as it implies $\lim_{t\to\infty}V(x\m^t,\lambda^t,\xi^t)=0$.
In view of Lemma~\ref{lem.V-Omega}, to prove the claim, it then suffices to show that $\tau\le T$. We prove this by contradiction. Assume $\tau> T$, then by definition of $\tau$ it holds
\begin{equation}\label{e.proof.uniform-conv.VT}
	V(x\m^T,\lambda^T,\xi^T)>\upsilon.
\end{equation}
As $(x^0,z^0,\lambda^0)\in\Xi_0$ implies $(x\m^0,\lambda^0,\xi^0)\in \Omega$,   Lemma~\ref{lem.V-Omega} implies $(x\m^t,\lambda^t,\xi^t)\in \Omega$ for all $t\in\N$. 
Moreover, Lemma~\ref{lem.V-Omega} and  $V(x\m^t,\lambda^t,\xi^t)>\upsilon$ for all $t< \tau$ imply
\begin{align*}
	\forall t\in\{0,\dots,\tau-1\},\quad V(x\m^{t+1},\lambda^{t+1},\xi^{t+1}) &\le V(x\m^t,\lambda^t,\xi^t) - \min\left\{ \frac{1}{8}\gamma\mu_f \varepsilon^2 ,\, \omega_\gamma \upsilon \right\}\\
	&= V(x\m^t,\lambda^t,\xi^t) - \omega_\gamma\upsilon,
\end{align*} 
where, in the last equality, we have used the definition of $\omega_\gamma$ in \eqref{d.omega_gamma} and  $\frac{1}{12}\gamma\mu_f \upsilon = \frac{1}{24}\gamma \mu_f \min\{\varepsilon^2,\eta^2\}\le \frac{1}{8}\gamma\mu_f\varepsilon^2$. 

As  $(x\m^0,\lambda^0,\xi^0)\in \Omega$ implies $V(x\m^0,\lambda^0,\xi^0)\le \frac{3}{2}\kappa_0^2$, under the assumption that 
 $\tau > T$, we then obtain
\begin{equation*}
	V(x\m^T,\lambda^T,\xi^T)\le \frac{3}{2}\kappa_0^2 - \omega_\gamma\upsilon T = \upsilon,
\end{equation*}
which contradicts \eqref{e.proof.uniform-conv.VT}.


\end{document}